\documentclass[letterpaper,11pt]{article}

\usepackage{ucs}
\usepackage[utf8x]{inputenc}
\usepackage{graphicx}
\usepackage{amsfonts}
\usepackage{dsfont}
\usepackage{amssymb}
\usepackage{amsmath}
\usepackage{amsthm}
\usepackage{enumerate}
\usepackage{stmaryrd}
\usepackage{fullpage}
\usepackage{ifthen}
\usepackage{subfigure}
\usepackage{epic}
\usepackage{authblk}
\usepackage{textcomp}
\usepackage{mathrsfs}
\usepackage{bm}
\usepackage[small]{caption}

\usepackage[hypertexnames=false,colorlinks=true,linkcolor=blue,citecolor=blue]{hyperref}
\usepackage[numbers,comma,square,sort&compress]{natbib}
\usepackage[letterpaper,text={6.5in,9in},centering]{geometry}

\usepackage{color}
\usepackage{titlesec}

\graphicspath{{eps/}{pdf/}}

\newcommand{\bqq}{\begin{equation}}
\newcommand{\eqq}{\end{equation}}
\newcommand{\bqs}{\begin{equation*}}
\newcommand{\eqs}{\end{equation*}}

\newcommand{\C}{\mathbb{C}} 
\newcommand{\D}{\mathbb{D}} 

\newcommand{\J}{\mathbb{J}}

\newcommand{\R}{\mathbb{R}} 

\newcommand{\W}{\mathbb{W}}

\renewcommand{\O}{\mathcal{O}}

\newcommand{\bG}{\mathbf{G}}
\newcommand{\bH}{\mathbf{H}}
\newcommand{\bK}{\mathbf{K}}
\newcommand{\cG}{\mathcal{G}}
\newcommand{\boldH}{\mathbf{H}}
\newcommand{\tcL}{\widetilde{\mathcal{L}}}

\newcommand{\cL}{\mathcal{L}}
\newcommand{\U}{\mathcal{U}}
\newcommand{\Q}{\mathcal{Q}}
\newcommand{\cW}{\mathcal{W}}
\newcommand{\cV}{\mathcal{V}}

\newcommand{\md}{\mathrm{d}}
\newcommand{\mbi}{\mathbf{i}}

\newcommand{\tp}{\tilde{p}}
\newcommand{\tq}{\tilde{q}}
\newcommand{\rmi}{\mathbf{i}}

\renewcommand{\bm}{\overline{m}(\lambda) }
\newcommand{\um}{\underline{m}(\lambda)}

\newcommand{\tscL}{\widetilde{\mathscr{L}}}
\newcommand{\scL}{\mathscr{L}}

\newcommand{\tscM}{\widetilde{\mathscr{M}}}
\newcommand{\scM}{\mathscr{M}}

\newtheorem{lem}{Lemma}[section]
\newtheorem{thm}{Theorem}
\newtheorem{prop}[lem]{Proposition}

\newenvironment{Proof}[1][.]%
 {\begin{trivlist}\item[]\textbf{Proof#1 }}%
 {\hspace*{\fill}$\rule{0.3\baselineskip}{0.35\baselineskip}$\end{trivlist}}

  {\begin{trivlist}\item[]{\bf Hypothesis #1 }\em}{\end{trivlist}}

\numberwithin{equation}{section}

\title{Asymptotic stability of the critical pulled front in a Lotka-Volterra competition model}

\author[1]{Gr\'egory Faye\footnote{email: \texttt{gregory.faye@math.univ-toulouse.fr}}}
\author[2]{Matt Holzer\footnote{email: \texttt{mholzer@gmu.edu}}}
\affil[1]{\small Institut de Math\'ematiques de Toulouse, UMR 5219, Universit\'e de Toulouse, UPS-IMT, F-31062 Toulouse Cedex 9 France}
\affil[2]{\small Department of Mathematical Sciences, George Mason University, Fairfax, VA 22030, USA}

\begin{document}
\maketitle

\begin{abstract}
We prove that the critical pulled front of Lotka-Volterra competition systems is nonlinearly asymptotically stable. More precisely, we show that perturbations of the critical front decay algebraically with rate $t^{-3/2}$ in a weighted $L^\infty$ space. Our proof relies on pointwise semigroup methods and utilizes in a crucial way that the faster decay rate $t^{-3/2}$ is a consequence of the lack of an embedded zero of the Evans function at the origin for the linearized problem around the critical front. 
\end{abstract}



\section{Introduction}

We consider the following system of reaction-diffusion equations known as the Lotka-Volterra competition model,
\bqq
\begin{split}
u_t&= u_{xx}+u(1-u-av),\\
v_t&=  \sigma v_{xx}+rv(1-bu-v),
\end{split}
\label{eq:main} 
\eqq
where $x\in\R$, $t>0$ and the non-negative parameters $a>0$, $b>0$, $\sigma>0$ and $r>0$ obey the following conditions, that we assume to be satisfied throughout the rest of the paper,
\bqq
0<a<1<b,
\label{par_mono}
\eqq
and 
\bqq
0<\sigma<2 \text{ and } (ab-M)r\leq M(2-\sigma)(1-a), \quad M:=\max\left\{1,2(1-a)\right\}.
\label{par_lin} 
\eqq
As a consequence of assumption \eqref{par_mono}, system (\ref{eq:main}) has three non-negative equilibrium points: the unstable state $(u,v)=(0,0)$, the unstable extinction state $(u,v)=(0,1)$ and the stable extinction state $(u,v)=(1,0)$. We are interested in non-negative traveling front solutions $(u(t,x),v(t,x))=(U_c(x-ct),V_c(x-ct))$ of \eqref{eq:main} connecting the two extinction states. Thus, it is follows that profiles $U_c(\xi)$ and $V_c(\xi)$ should satisfy
\bqq
\begin{split}
U''+cU'+U(1-U-aV)&=0,\\
\sigma V''+cV'+rV(1-bU-V)&=0,
\end{split}
\label{eq:wave} 
\eqq
subject to
\bqq
\underset{\xi\rightarrow-\infty}{\lim}(U_c(\xi),V_c(\xi))=(1,0), \text{ and } \underset{\xi\rightarrow+\infty}{\lim}(U_c(\xi),V_c(\xi))=(0,1).
\label{eq_limits} 
\eqq

The existence of such fronts has been studied extensively \cite{kanon97,lewis02,weinberger02,fei03,huang10,alhasanat18,morita09}. In particular, it is known that there exists $c_{\text{min}}\geq 0$ so that non-negative traveling front solutions of \eqref{eq:wave}-\eqref{eq_limits} exist if and only if $c\geq c_{\text{min}}$, see \cite{kanon97,lewis02,weinberger02}. We shall refer to $c_{\text{min}}$ as the minimal wavespeed and traveling fronts propagating at the minimal speed are called critical. Due to the invariance by translation of \eqref{eq:main}, there is a one-parameter family of such traveling fronts for each $c\geq c_{\text{min}}$.  A direct linearization around the unstable extinction state $(u,v)=(0,1)$ reveals that the minimal wavespeed  satisfies
\bqs 
c_{\text{min}} \geq c_*:= 2\sqrt{1-a}. 
\eqs
We call $c_*=2\sqrt{1-a}$ the linear spreading speed, as it is the asymptotic spreading speed of compactly supported initial data for the system linearized near the unstable state $(0,1)$. It was shown in \cite{lewis02,weinberger02,alhasanat18} that for any parameters satisfying \eqref{par_mono}-(\ref{par_lin}), $c_{\text{min}}=c_*$ such that the minimal speed is linearly selected. We denote by $(U_*,V_*)$ the corresponding critical front profile solution of \eqref{eq:wave}-\eqref{eq_limits} propagating with the linear spreading speed $c_*=2\sqrt{1-a}$. We note that the existence or non-existence of a critical front has been the focus of a great deal of research.  We refer the reader to \cite{fei03,huang10} for other works on the existence of critical traveling fronts propagating with the linear speed $c_*$.  It bears mentioning that there are choices of parameters for which the critical front propagates with minimal speed $c_{\text{min}}$ strictly larger than $c_*$, see for example \cite{hosono03,huang11,holzer12}. In the following, in order to avoid any confusion with the later case of nonlinear determinacy $c_{\text{min}}>c_*$, we shall use the terminology popularized in \cite{vansaarloos03} and refer to the critical front $(U_*,V_*)$ propagating with the linear speed $c_*=2\sqrt{1-a}$ as the critical pulled front.

The stability of traveling front solutions propagating into unstable states has been the subject of numerous studies in the past decades.  In the scalar case, the stability of super-critical fronts was established in \cite{sattinger} using exponential weights to stabilize the essential spectrum.  For the critical pulled front, stability was first established by \cite{kirchgassner92}, and later extended and refined in \cite{bricmont92,eckmann94,gallay94}. The sharpest of these results proves that perturbations of the critical pulled front converge in an exponentially weighted $L^\infty$ space with algebraic rate $t^{-3/2}$, see \cite{gallay94,faye19}. Of course, using comparison principle techniques, strong results concerning the convergence of compactly supported initial data to traveling fronts are possible; see for example \cite{lewis02}. 

The literature concerning the stability of pulled fronts in systems of reaction-diffusion equations is smaller.  We mention the results of \cite{gardner82,kanon96} for traveling fronts connecting stable extinction states. For the Lotka-Volterra system (\ref{eq:main}) under the assumptions \eqref{par_mono}, the local stability of traveling fronts propagating with wavespeed $c>c_{\text{min}}$ for \eqref{eq:main} was proved in \cite{leung11}.  The strategy employed there is to use exponential weights to stabilize the essential spectrum and comparison principles to exclude unstable point spectrum.  Since $c>c_{min}$ these fronts have weak exponential decay and so the derivative of the front does not lead to a translational eigenvalue in the weighted space and the spectrum lies strictly in the negative half plane.  Taken together, this leads to exponential in time stability with respect to spatially localized perturbations.

In the present work, we are concerned with the asymptotic stability of the critical pulled front $(U_*,V_*)$. More precisely, we will show that perturbations of the critical pulled front converge in an exponentially weighted $L^\infty$ space with algebraic rate $t^{-3/2}$, as it is the case for the scalar Fisher-KPP equation \cite{gallay94,faye19}. It is important to emphasize that the essential spectrum of the linear operator obtained by linearizing \eqref{eq:main} in a coordinate frame moving at speed $c_*$ cannot be stabilized by using exponential weights. In fact, using appropriate exponential weights, one finds continuous spectrum up to the origin.  Essential to the analysis here is the absence of an embedded eigenvalue at the origin.  The lack of an embedded translational eigenvalue at the origin is expected due to the fact that the critical pulled front has weak exponential decay near $\xi=+\infty$,
\bqq \left( \begin{array}{c} U_*(\xi) \\ V_*(\xi)-1 \end{array}\right) \sim \xi e^{-\gamma_* \xi}\left( \begin{array}{c} \beta_1 \\ \beta_2 \end{array}\right) +o\left(\xi e^{-\gamma_* \xi}\right) \quad \text{as $\xi\to\infty$}, 
\label{eq:asympt_front}
\eqq
where $\gamma_*:=c_*/2=\sqrt{1-a}$ and for some vector $(\beta_1,\beta_2)\neq 0$, see \cite{kanon97,morita09,girardin18}. As already noticed in \cite{sandstede04} and further exploited in \cite{faye19}, this weak exponential decay implies that the derivative of the critical pulled front has also weak exponential decay and thus does not contribute to a zero eigenvalue of the linearized operator.

We now state our main result. First, we rewrite system \eqref{eq:main} in traveling wave coordinate  moving to the right with speed $c_*$
\bqq
\begin{split}
u_t&= u_{xx}+c_*u_x+u(1-u-av),\\
v_t&=  \sigma v_{xx}+c_*v_x+rv(1-bu-v),
\end{split}
\label{eq:mainTW} 
\eqq
where $x\in\R$, $t>0$, such that the critical pulled front $(U_*,V_*)$ is now a stationary solution \eqref{eq:mainTW}. We introduce a positive, bounded, smooth weight function $\omega(x)>0$ of the form
\bqq
\omega(x)=\left\{ 
\begin{array}{ll}
e^{-\gamma_*x} & x \geq 1,\\
e^{\delta x} & x \leq -1,
\end{array} \right. 
\label{eq:weights}
\eqq
for some $\delta>0$ that satisfies
\bqq
0<\delta< \min\left(-\gamma_*+\sqrt{\gamma_*^2+1},\frac{-\gamma_*+\sqrt{\gamma_*^2+\sigma r(b-1)}}{\sigma}\right).
\label{eq:defdb}
\eqq

\begin{thm}\label{thmmain}
Assume that assumptions \eqref{par_mono}-\eqref{par_lin} hold. Consider \eqref{eq:mainTW} with initial datum $(u(0,\cdot),v(0,\cdot))=(U_*,V_*)+(u_0,v_0)$ satisfying 
\bqs
0\leq u(0,x) \leq 1, \text{ and }  0\leq v(0,x) \leq 1, \quad \forall x\in\R.
\eqs
There exists $C>0$ and $\epsilon>0$ such that if $(p_0,q_0):=\left(u_0/\omega,v_0/\omega\right)$ satisfies 
\bqs
\left\| (p_0,q_0)  \right\|_{L^\infty(\R)}+\left\| (1+|\cdot|) (p_0,q_0)  \right\|_{L^1(\R)} <\epsilon,
\eqs
then the solution $(u(t,x),v(t,x))$ is defined for all time and the critical pulled front is nonlinearly stable in the sense that
\bqs
\left\| \frac{1}{(1+|\cdot|)} \left( p(t,\cdot),q(t,\cdot) \right)  \right\|_{L^\infty(\R)} \leq \frac{C\epsilon}{(1+t)^{3/2}}, \quad t>0,
\eqs
where $p(t,x):=(u(t,x)-U_*(x))/\omega(x)$ and $q(t,x):=(v(t,x)-V_*(x))/\omega(x)$ for all $x\in\R$ and $t>0$.
\end{thm}

Our strategy of proof follows the lines of the one we recently proposed for the scalar Fisher-KPP equation \cite{faye19}. It is based upon pointwise semigroup methods which were introduced in \cite{zumbrun98}, and have been  developed precisely to address stability problems where the essential spectrum cannot be separated from the imaginary axis. We refer to \cite{zumbrun98,howard02,johnson11,beck14} for various applications of pointwise semigroup method to the stability of viscous shock waves, stability and instability of spatially periodic patterns or  stability of defects in reaction-diffusion equations, to mention a few.   Most related to the present study is \cite{howard07} where the stability of the kink solution for the Cahn-Hilliard equation is studied.  In a nutshell, our approach can be summarized as follows.
\begin{itemize}
\item We construct bounded linearly independent solutions $(\phi_{1}^-,\phi_{2}^-)$ and $(\phi_{1}^+,\phi_{2}^+)$ to the eigenvalue problem 
\[ \mathcal{L} \left(\begin{matrix} p \\ q \end{matrix}\right)=\lambda \left(\begin{matrix} p\\ q \end{matrix}\right),\]
on $ \R^-$ and $\R^+$ respectively,  where $\mathcal{L}$ is a linear operator describing the linearized eigenvalue problem of \eqref{eq:mainTW} around $(U_*,V_*)$ transformed to a weighted space with weight given by $\omega$ from \eqref{eq:weights}.
\item We find bounds on the pointwise Green's function $\bG_\lambda(x,y)$ which is a $2\times2$ matrix satisfying 
\bqs \left(\mathcal{L}-\lambda \mathrm{I}_2\right)  \bG_\lambda(x,y)=-\delta(x-y)\mathrm{I}_2, \quad x,y \in\R,\eqs
and that can be expressed in terms of combinations of $\phi_{1,2}^\pm$ evaluated at $x$ and/or $y$.
\item We apply the inverse Laplace transform, and show that temporal Green's function
\bqs
\cG(t,x,y)=\frac{1}{2\pi \mathbf{i}} \int_\mathscr{C} e^{\lambda t}\bG_\lambda(x,y)
\mathrm{d}\lambda
\eqs
decays pointwise with algebraic rate $t^{-3/2}$ by a suitable choice of inversion contour $\mathscr{C} \subset \mathbb{C}$.
\item The final step consists in applying $L^p$ estimates to the nonlinear solution expressed using Duhamel's formula 
\bqs
\left(\begin{matrix}p(t,x) \\ q(t,x) \end{matrix}\right)  = \int_\R \cG(t,x,y)\left(\begin{matrix}p_0(y) \\ q_0(y) \end{matrix}\right)\mathrm{d}y+\int_0^t \int_\R \cG(t-\tau,x,y)\left( \begin{matrix}  \mathcal{N}_u(p(\tau,y),q(\tau,y)) \\  \mathcal{N}_v (p(\tau,y),q(\tau,y))  \end{matrix}\right) \mathrm{d}y\mathrm{d}\tau,
\eqs
where $(\mathcal{N}_u(p,q),\mathcal{N}_v(p,q))$ encodes nonlinear terms, and show that the nonlinear system also exhibits the same algebraic decay rate.
\end{itemize}

Of the four steps above, we expend the most effort obtaining the necessary bounds on the pointwise Green's function $\bG_\lambda(x,y)$.   This involves two pieces.  First, we must show that $\bG_\lambda(x,y)$ has no poles for $\lambda$ in the right half plane.  These eigenvalues would necessarily imply instability of the front.  To do this, we exploit the monotone structure of system (\ref{eq:main}) and apply comparison principle techniques to rule out unstable point spectrum; see \cite{bates2006spectral,leung11}.  Next, we require that $\bG_\lambda(x,y)$ can be expressed as  $e^{-\sqrt{\lambda}|x-y|}\mathbf{H}_\lambda(x,y)$, with $\mathbf{H}_\lambda(x,y)$ a $2\times2$ matrix, analytic in $\lambda$ in the region considered and bounded uniformly as a function of $x$ and $y$.  The prefactor $e^{-\sqrt{\lambda}|x-y|}$ is the Laplace transform of the derivative of the heat kernel in one space dimension; from which the decay rate $t^{-3/2}$ is naturally expected.  The lack of a singularity in $\bG_\lambda(x,y)$ as $\lambda\to 0$ is enforced by the lack of an embedded eigenvalue at the edge of the essential spectrum.  This, in turn, is related to the weak exponential decay of the critical front which precludes the derivative of the front from contributing a singularity at $\lambda=0$ and the monotone structure of (\ref{eq:main}) which can be used to preclude any other bounded solution to the eigenvalue problem at $\lambda=0$.  Once these bounds have been established the third and fourth steps in our proof are closely related to the scalar case presented in \cite{faye19}.

We expect that the method used to establish Theorem~\ref{thmmain} could be used to establish nonlinear stability of other pulled fronts.  In particular, we anticipate that the two main ingredients necessary for the proof: the lack of an embedded translational eigenvalue at $\lambda=0$ and a subsequent bound on the pointwise Green's function of the form $e^{-\sqrt{\lambda}|x-y|}$ would occur generically in examples of pulled fronts.  However, presenting a significant challenge from a mathematical point of view is excluding unstable point spectrum.  The principle tool in this case is the Evans function; we refer the reader to \cite{kapitula,sandstede02} for an introduction.  

The rest of the paper is organized as follows. In Section 2, we set up and study the linearized
eigenvalue problem. In Section 3, we derive bounds on the pointwise Green’s function $\bG_\lambda(x,y)$.  In Section 4, we use these estimates to obtain estimates for the temporal Green's function $\cG(t,x,y)$.  In Section 5, we use material from the previous sections to prove nonlinear stability.  Some proofs are given in the Appendix.

\section{Preliminaries and ODE estimates}\label{secprelim}

Consider \eqref{eq:mainTW} and write the solutions
\[ \left(\begin{array}{c} u(t,x) \\ v(t,x) \end{array}\right)= \left(\begin{array}{c} U_*(x) \\ V_*(x) \end{array}\right)+\left(\begin{array}{c} \tp(t,x) \\ \tq(t,x) \end{array}\right), \quad t>0 \text{ and } x\in\R. \]
Using the fact that the profile $(U_*,V_*)$ of the critical pulled front verifies \eqref{eq:wave}, we derive the following set of equations for the perturbations 
\bqq
\begin{split}
\tp_t&= \tp_{xx}+c_*\tp_x + \tp(1-2U_*(x)-aV_*(x))-aU_*(x)\tq -\tp(\tp+a\tq), \\
\tq_t&=  \sigma \tq_{xx}+c_*\tq_x  -rbV_*(x)\tp +r\tq(1-bU_*(x)-2V_*(x)) -r\tq(b\tp+\tq),
\end{split}
 \label{eq:linfirst}
\eqq
for $t>0$ and $x\in\R$. This system can be expressed as
\[ \partial_t \left(\begin{array}{c} \tp \\ \tq \end{array}\right)= \widetilde{\mathcal{L}} \left(\begin{array}{c} \tp \\ \tq \end{array}\right) +\widetilde{\mathcal{N}}(\tp,\tq), \quad t>0 \text{ and } x\in\R,\]
with
\bqs
\widetilde{\mathcal{L}}:=\left( \begin{array}{cc} \partial_{xx}+c_*\partial_x+(1-2U_*(x)-aV_*(x)) & -aU_*(x) \\
-rbV_*(x) & \sigma\partial_{xx}+c_*\partial_x+r(1-bU_*(x)-2V_*(x)) \end{array}\right),
\eqs
and
\bqs
\widetilde{\mathcal{N}}(\tp,\tq):=\left( \begin{array}{c}  -\tp(\tp+a\tq) \\  -r\tq(b\tp+\tq) \end{array}\right).
\eqs
As the profile $(U_*,V_*)$ satisfies the asymptotics \eqref{eq_limits}, the linear operator $\widetilde{\mathcal{L}}$ is exponentially asymptotic as $x\to\pm\infty$ to the following limiting operators
\begin{align*}
 \tilde{\mathcal{L}}_+&:=\left( \begin{array}{cc} \partial_{xx}+c_*\partial_x+(1-a) & 0 \\
-rb & \sigma\partial_{xx}+c_*\partial_x-r \end{array}\right), \\
 \tilde{\mathcal{L}}_-&:=\left( \begin{array}{cc} \partial_{xx}+c_*\partial_x-1 & -a \\
0 & \sigma\partial_{xx}+c_*\partial_x+r(1-b) \end{array}\right).
\end{align*}
The essential spectrum of $\widetilde{\mathcal{L}}_+$ is unstable (in $L^2(\mathbb{R})$ for instance) since $1-a>0$ in the upper left component as a consequence of assumption \eqref{par_mono}.  The essential spectrum can be shifted using exponential weights and to this end we introduce the weight $\omega$, defined in \eqref{eq:weights}, and that writes
\bqs
\omega(x)=\left\{ 
\begin{array}{ll}
e^{-\gamma_*x} & x \geq 1,\\
e^{\delta x} & x \leq -1,
\end{array} \right. 
\eqs
for some well chosen $\delta>0$ and $\gamma_*=c_*/2=\sqrt{1-a}>0$. The conditions on $\delta$, as prescribed by \eqref{eq:defdb}, will be explained below. Without loss of generality we further assume that $\omega(0)=1$.   

Let 
\[ \tp(t,x)= \omega(x) p(t,x), \quad \tq(t,x)=\omega(x)q(t,x), \quad t>0 \text{ and } x\in\R. \]
This converts (\ref{eq:linfirst}) into
\bqq
\begin{split}
p_t=&~ p_{xx}+\left(c_*+2\frac{\omega'}{\omega}\right) p_x +p\left(1+c_*\frac{\omega'}{\omega}+\frac{\omega''}{\omega}-2U_*(x)-aV_*(x)\right) -a U_*(x) q\\
& - \omega p( p+aq),   \\
q_t=&~  \sigma q_{xx}+\left(c_*+2\sigma \frac{\omega'}{\omega}\right) q_x
+rq\left(1 +c_*\frac{\omega'}{r\omega}+\frac{\sigma\omega''}{r\omega}-bU_*(x)-2V_*(x)\right) \\
 &  -rbV_*(x) p -r\omega q(bp+q), 
\end{split}
\label{eq:pq} 
\eqq
for $t>0$ and $x\in\R$.
We will find it convenient to introduce the notation,
\[ \partial_t \left(\begin{array}{c} p \\ q \end{array}\right)= \mathcal{L} \left(\begin{array}{c} p\\ q \end{array}\right) +\mathcal{N}(p,q), \quad t>0,\quad  x\in\R, \]
where
\bqs
\mathcal{L}:=\left( \begin{array}{cc} \mathcal{L}_u & \cL_{12} \\
\cL_{21} & \mathcal{L}_v \end{array}\right), \text{ and } \mathcal{N}(p,q):=\omega \left( \begin{array}{c}  - p(p+a q) \\  -rq(b p+ q)\end{array}\right),
\eqs
with
\begin{align*}
\mathcal{L}_u&:=\partial_{xx}+\left(c_*+2\frac{\omega'}{\omega}\right) \partial_x +\left(1+c_*\frac{\omega'}{\omega}+\frac{\omega''}{\omega}-2U_*(x)-aV_*(x)\right),\\
\cL_{12}&:=-aU_*(x),\\
\mathcal{L}_v&:=\sigma \partial_{xx}+\left(c_*+2\sigma \frac{\omega'}{\omega}\right) \partial_x
+r\left(1 +c_*\frac{\omega'}{r\omega}+\frac{\sigma\omega''}{r\omega}-bU_*(x)-2V_*(x)\right),\\
\cL_{21}&:=-rbV_*(x).
\end{align*}
For $x\geq 1$, system (\ref{eq:pq})  reduces to
\bqq
\begin{split}
p_t=&~ p_{xx} -p(2U_*(x)+a(V_*(x)-1)) -aU_*(x) q - p( p+aq)e^{-\gamma_* x},  \\
q_t=&~  \sigma q_{xx}+c_*\left(1-\sigma\right) q_x
+rq\left(-1+(\sigma-2)\frac{\gamma_*^2}{r}-bU_*(x)-2(V_*(x)-1)\right)  \\
 & -rbV_*(x) p -rq(bp+q)e^{-\gamma_* x}. 
\end{split} 
 \label{eq:pq+}
\eqq
For $x\leq -1$, system (\ref{eq:pq})  reduces to
\bqq
\begin{split}
p_t=&~ p_{xx} +\left(c_*+2\delta \right) p_x+p(-1+c_*\delta+\delta^2-2(U_*(x)-1)-aV_*(x)) -aU_*(x) q  \\
& -p(  p+a q)e^{\delta x},  \\
q_t=&~  \sigma q_{xx}+\left(c_*+2\sigma\delta \right) q_x
+rq\left(1-b +\frac{c_*\delta}{r} +\sigma \frac{\delta^2}{r}-b(U_*(x)-1)-2V_*(x)\right)- rbV_*(x) p \\
  & -rq(bp+ q)e^{\delta x}. 
\end{split}
\label{eq:pq-} 
\eqq 
Here, we need to select $\delta>0$ sufficiently small to ensure that both essential spectra of $\mathcal{L}_u$ and $\mathcal{L}_v$ remains stable. Thus, we both need that
\bqs
-1+c_*\delta+\delta^2<0 \text{ and } 1-b +\frac{c_*\delta}{r} +\sigma \frac{\delta^2}{r}<0,
\eqs
and we select $\delta$ such that
\bqs
0<\delta < \min\left(-\gamma_*+\sqrt{\gamma_*^2+1},\frac{-\gamma_*+\sqrt{\gamma_*^2+\sigma r(b-1)}}{\sigma}\right),
\eqs
which is precisely the condition in \eqref{eq:defdb}.

As a consequence of the above analysis, we have that
\bqs
\cL:\mathcal{D}(\cL)\subset L^2(\R)\times L^2(\R)\rightarrow L^2(\R)\times L^2(\R),
\eqs
with dense domain $\mathcal{D}(\cL)=H^2(\R)\times H^2(\R)$ is a closed operator with Fredholm borders, denoted $\sigma_F(\cL)$, that are those curves in the complex $\lambda$-plane defined as
\bqs
\sigma_F(\cL):=\left\{\lambda \in \mathbb{C} ~|~ d_\pm(\mbi \ell, \lambda)=0, \text{ for some } \ell \in\R \right\},
\eqs
with 
\begin{align*}
d_-(\nu,\lambda)&:=\left(\nu^2 +\left(c_*+2\delta \right) \nu-1+c_*\delta+\delta^2-\lambda \right)\left(\sigma \nu^2+\left(c_*+2\sigma\delta \right) \nu
+r(1-b) +c_*\delta +\sigma \delta^2-\lambda\right),\\
d_+(\nu,\lambda)&:=\left(\nu^2-\lambda \right)\left(\sigma \nu^2+c_*\left(1-\sigma\right) \nu
-r +(\sigma-2) \gamma_*^2-\lambda\right).
\end{align*}
Here, the Fredholm borders of $\cL$ are the union of four curves in the complex plane
\bqs
\sigma_F(\cL)=\Gamma^-_u\cup\Gamma^-_v\cup\Gamma^+_u\cup\Gamma^+_v,
\eqs
given by
\begin{align*}
\Gamma^-_u&:=\left\{-\ell^2-1+c_*\delta+\delta^2+\left(c_*+2\delta \right)\mbi\ell~|~ \ell \in\R\right\} ,\\
\Gamma^-_v&:=\left\{-\sigma\ell^2+r(1-b) +c_*\delta +\sigma \delta^2+\left(c_*+2\sigma\delta \right)\mbi\ell ~|~ \ell \in\R\right\},\\
\Gamma^+_u&:=]-\infty,0],\\
\Gamma^+_v&:=\left\{-\sigma\ell^2-r +(\sigma-2) \gamma_*^2+c_*\left(1-\sigma \right) \mbi\ell~|~ \ell \in\R \right\}.
\end{align*}
Our careful choice of $\delta$ ensures that the three parabola $\Gamma^-_u$, $\Gamma^-_v$ and $\Gamma^+_v$ are entirely located to the left of the imaginary axis in the complex plane. Only the half-line $\Gamma^+_u$ touches the imaginary axis at the origin. Finally we set $\iota:=\max\left\{-1+c_*\delta+\delta^2,r(1-b) +c_*\delta +\sigma \delta^2,-r +(\sigma-2) \gamma_*^2 \right\}<0$ and introduce the vertical line in the complex plane:
\bqs
\Gamma:=\left\{ \iota + \mbi\ell ~|~ \ell \in\R\right\}.
\eqs

\subsection{The linearized eigenvalue problem}
Our eventual goal is to derive bounds on the temporal Green's function $\cG(t,x,y)$ via bounds on the pointwise Green's function $\bG_\lambda(x,y)$.  The pointwise Green's function is constructed from exponentially decaying solutions to the linearized eigenvalue problem,
\[ \mathcal{L} \left(\begin{array}{c} p\\ q \end{array}\right)=\lambda \left(\begin{array}{c} p\\ q \end{array}\right), \quad \lambda \in\C, \quad p,q \in H^2(\R).\]
This system is equivalent to a system of two coupled second order ODEs.  To study its solutions, we express it as a system of four first order equations and introduce the notation
\bqq
P'=\mathcal{A}(x,\lambda)P, \quad P=(p,p',q,q')^{\mathbf{t}},
\label{eqPsyst}
\eqq
where 
\bqs
\mathcal{A}(x,\lambda):=\left(\begin{matrix} 0 & 1 & 0 & 0 \\
 \lambda-\zeta_u(x) & -\left(c_*+2\dfrac{\omega'(x)}{\omega(x)}\right) & aU_*(x)& 0 \\
0 & 0 & 0 & 1 \\
\dfrac{rb}{\sigma}V_*(x)& 0 & \dfrac{\lambda-\zeta_v(x)}{\sigma} & -\dfrac{1}{\sigma}\left(c_*+2\sigma \dfrac{\omega'(x)}{\omega(x)}\right) \end{matrix} \right),
\eqs
where
\begin{align*}
\zeta_u(x)&:=1+c_*\dfrac{\omega'(x)}{\omega(x)}+\dfrac{\omega''(x)}{\omega(x)}-2U_*(x)-aV_*(x),\\
\zeta_v(x)&:=r\left(1 +c_*\frac{\omega'(x)}{r\omega(x)}+\frac{\sigma\omega''(x)}{r\omega(x)}-bU_*(x)-2V_*(x)\right).
\end{align*}
It is important to note that the matrix $\mathcal{A}(x,\lambda)$ simplifies for $x\geq1$ and $x\leq-1$ according to \eqref{eq:pq+} and \eqref{eq:pq-}. More precisely, we have that
\bqq
\mathcal{A}(x,\lambda)=\left\{\begin{array}{cc}
\mathcal{A}^+(\lambda)+\mathcal{B}^+(x), & x\geq 1, \\
\mathcal{A}^-(\lambda)+\mathcal{B}^-(x), & x\leq -1,
\end{array}\right. \label{eq:Aeq}
\eqq
with associated asymptotic matrices at $x=\pm\infty$
\begin{align*}
\mathcal{A}^+(\lambda)&:=\left(\begin{matrix} 0 & 1 & 0 & 0 \\ 
\lambda & 0  & 0 & 0 \\ 
0 & 0  & 0 & 1 \\ 
\frac{rb}{\sigma} & 0 & \dfrac{\lambda+r+(2-\sigma)\gamma_*^2}{\sigma} & \dfrac{c_*(\sigma-1) }{\sigma}  \end{matrix} \right), \\
\mathcal{A}^-(\lambda)&:=\left(\begin{matrix} 0 & 1 & 0 & 0 \\ 
\lambda-(-1+c_*\delta+\delta^2) & -(c_*+2\delta)  & a & 0 \\ 0 & 0  & 0 & 1 \\ 
0 & 0 & \dfrac{\lambda -(r(1-b)+c_*\delta+\sigma\delta^2)}{\sigma} & - \dfrac{c_*+2\sigma \delta }{\sigma} \end{matrix} \right),
\end{align*}
and associated remainder matrices given by
\begin{align*}
\mathcal{B}^+(x)&:=\left(\begin{matrix} 0 & 0 & 0 & 0 \\ 
2U_*(x)+a(V_*(x)-1) & 0  & aU_*(x) & 0 \\ 
0 & 0  & 0 & 0 \\ 
\dfrac{rb}{\sigma}(V_*(x)-1) & 0 & \dfrac{2r(V_*(x)-1)+brU_*(x)}{\sigma} & 0 \end{matrix} \right),\\
\mathcal{B}^-(x)&:=\left(\begin{matrix} 0 & 0 & 0 & 0 \\ 
2(U_*(x)-1)+aV_*(x) & 0  & a(U_*(x)-1) & 0 \\ 
0 & 0  & 0 & 0 \\ 
\dfrac{rb}{\sigma}V_*(x) & 0 & \dfrac{2rV_*(x)+rb(U_*(x)-1)}{\sigma} & 0 \end{matrix} \right).
\end{align*}

Furthermore, there exists $0<\alpha<\min\{\delta,\gamma_*\}$ such that
\bqs
\| \mathcal{B}^+(x) \| \leq C_+ e^{-\alpha x}, \quad x\geq 1,\quad \text{ and }\quad \| \mathcal{B}^-(x) \| \leq C_- e^{\alpha x}, \quad x\leq -1,
\eqs
for some constants $C_\pm>0$.

The matrix $\mathcal{A}^+(\lambda)$ has four eigenvalues
\[ \pm\sqrt{\lambda}, \quad \nu_v^\pm(\lambda)=\gamma_* -\frac{\gamma_*}{\sigma}\pm\frac{1}{\sigma}\sqrt{\gamma_*^2+\sigma (r +\lambda)}, \]
and four eigenvectors
\[ e_u^\pm(\lambda) =\left(\begin{array}{c} 1 \\ \pm \sqrt{\lambda} \\ y_v^\pm(\lambda) \\ \pm\sqrt{\lambda} y_v^\pm(\lambda) \end{array}\right), \quad e_v^\pm(\lambda) =\left(\begin{array}{c} 0 \\ 0 \\ 1 \\ \nu_v^\pm(\lambda) \end{array}\right), \]
with
\bqq
y_v^\pm(\lambda)=\frac{rb}{ \sigma\lambda\pm\sqrt{\lambda} c_*(1-\sigma)-\lambda-r-(2-\sigma)\gamma_*^2}. \label{eq:yv}
\eqq
We note that $y_v^\pm(\lambda)$ is bounded so long as $\nu_v^\pm(\lambda)\neq \pm \sqrt{\lambda}$.  In the following Lemma, we characterize solutions of system \eqref{eqPsyst} with prescribed exponential growth and decay rates at $+\infty$. We restrict $\lambda \in \C$ to a subset of the complex plane near the origin, to the right of $\Gamma$ and off the negative real axis, for which there exist positive constants $\eta_\pm>0$ such that 
\bqq
\mathrm{Re}(\nu_v^-(\lambda)) <-\eta_- < -\mathrm{Re}(\sqrt{\lambda})  \leq 0 \leq \mathrm{Re}(\sqrt{\lambda})<\eta_+ < \mathrm{Re}(\nu_v^-(\lambda)). 
\label{eq:orderingeig}
\eqq
 Note that without loss of generality we can assume $-\eta_-<-\alpha<0$ since $\eta_-$ is strictly positive for the $\lambda$ values considered here and there is no harm in taking $\alpha$ smaller. Indeed, there always exists $M_s>0$, small enough, such that for all $|\lambda|<M_s$, to the right of $\Gamma$ and off the negative real axis, the ordering \eqref{eq:orderingeig} of the eigenvalues of $\mathcal{A}^+(\lambda)$ is valid.

\begin{lem} \label{lem:+solsp} For $x\geq 0$ and for all $|\lambda|<M_s$, to the right of $\Gamma$ and off the negative real axis, we have the existence of two linearly independent bounded solutions $\phi_{1,2}^+$ together with two linearly independent unbounded solutions $\psi_{1,2}^+$ with prescribed asymptotic growth rates $\sqrt{\lambda}$ and $\nu_v^+(\lambda)$, 
\bqq
\begin{split}
\phi_1^+(x)&= e^{-\sqrt{\lambda} x}\left( e_u^-(\lambda)+\theta_1^+(x,\lambda)\right) \\
 \phi_2^+(x)&= e^{\nu_v^-(\lambda) x}\left( e_v^-(\lambda)+\theta_2^+(x,\lambda)\right) \\
\psi_1^+(x)&= e^{\sqrt{\lambda} x}\left( e_u^+(\lambda)+\kappa_1^+(x,\lambda)\right) \\
\psi_2^+(x)&= e^{\nu_v^+(\lambda) x}\left( e_v^+(\lambda)+\kappa_2^+(x,\lambda)\right)
\end{split}
\label{eq:+solsp}
\eqq
Furthermore, we have the following uniform bounds for the the vectors in (\ref{eq:+solsp}),
\bqs
\left\|\theta_{1,2}^+(x,\lambda)\right\|\leq C e^{-\alpha x} \text{ and } \left\|\kappa_{1,2}^+(x,\lambda)\right\|\leq C e^{-\alpha x}, \quad x\geq0, 
\eqs 
 for some constant $C>0$ independent of $\lambda$.
\end{lem}

\begin{Proof} This result follows from Proposition 3.1 of \cite{zumbrun98}.  For the convenience of the reader and because we will require some details of the proof to establish a subsequent lemma we include a discussion of the proof here.  We will demonstrate how to construct $\phi_1^+(x)$ and $\psi_1^+(x)$.  Associated to each eigenvalue of the matrix $\mathcal{A}^+(\lambda)$ is a spectral projection which we denote by $\Pi_{ss}(\lambda)$, $\Pi_{ws}(\lambda)$, $\Pi_{wu}(\lambda)$, $\Pi_{su}(\lambda)$ for the eigenvalues $\nu_v^-(\lambda)$, $-\sqrt{\lambda}$, $\sqrt{\lambda}$, $\nu_v^+(\lambda)$ respectively. We also use the notation $\Pi_{cu}(\lambda)=\mathrm{Id}-\Pi_{ss}(\lambda)$ for the center-unstable projection.

We first construct $\phi_1^+(x)$.  First, we transform  (\ref{eqPsyst}) by letting $Z(x)=e^{-\sqrt{\lambda}x}P(x)$ and note that $Z(x)$ obeys
\bqq
 Z'(x)=\left(\mathcal{A}^+(\lambda)+\sqrt{\lambda}I_4\right)Z(x)+ \mathcal{B}^+(x)Z(x),\quad x\geq1.
 \label{eqZphi}
 \eqq
The solution of this equation can be expressed as 
\[ Z(x)=e^{(\mathcal{A}^+(\lambda)+\sqrt{\lambda}I_4)x}Z(x_0)+\int_{x_0}^x e^{(\mathcal{A}^+(\lambda)+\sqrt{\lambda}I_4)(x-y)} \mathcal{B}^+(y)Z(y)\md y, \quad x\geq1, \]
for some $x_0\geq1$. We take 
\[ Z(x_0)=e_u^-(\lambda)-\int_{x_0}^\infty e^{-(\mathcal{A}^+(\lambda)+\sqrt{\lambda}I_4)y}\Pi_{cu}(\lambda) \mathcal{B}^+(y)Z(y)\md y, \]
such that 
\bqq 
\begin{split}
Z(x)=&~e_u^-(\lambda)+\int_{x_0}^x e^{(\mathcal{A}^+(\lambda)+\sqrt{\lambda}I_4)(x-y)}\Pi_{ss}(\lambda) \mathcal{B}^+(y)Z(y)\md y\\
& - \int_{x}^\infty e^{(\mathcal{A}^+(\lambda)+\sqrt{\lambda}I_4)(x-y)}\Pi_{cu}(\lambda) \mathcal{B}^+(y)Z(y)\md y.
\end{split}
\label{eqZphimap+} \eqq
There exists a constant $K>0$, independent of $\lambda$ for $|\lambda|<M_s$, and to the right of $\Gamma$ such that the strong stable projection $\Pi_{ss}(\lambda)$ and the center-unstable projection $\Pi_{cu}(\lambda)$ obey the following bounds
\begin{eqnarray*} \left\|e^{(\mathcal{A}^+(\lambda)+\sqrt{\lambda}I_4)x}\Pi_{ss}(\lambda)\right\|&\leq& K e^{(-\eta_-+\mathrm{Re}(\sqrt{\lambda}))x}, \quad x>0, \\
\left\|e^{(\mathcal{A}^+(\lambda)+\sqrt{\lambda}I_4)x}\Pi_{cu}(\lambda)\right\|&\leq& K,  \quad x<0,
\end{eqnarray*} 
which when taken together with $\left\|\mathcal{B}^+(x)\right\|\leq C_+ e^{-\alpha x}$ implies that the mapping (\ref{eqZphimap+}) is a contraction mapping on $L^{\infty}([x_0,\infty))$ for some $x_0$ sufficiently large.  We can then expand $Z(x)=e_u^-(\lambda)+\theta_1(x,\lambda)$ with 
\bqq
\begin{split} \theta_1^+(x,\lambda)=& \int_{x_0}^x e^{(\mathcal{A}^+(\lambda)+\sqrt{\lambda}I_4)(x-y)}\Pi_{ss}(\lambda) \mathcal{B}^+(y)(e_u^-(\lambda)+\theta_1^+(y,\lambda))\md y  \\
&- \int_{x}^\infty e^{(\mathcal{A}^+(\lambda)+\sqrt{\lambda}I_4)(x-y)}\Pi_{cu}(\lambda) \mathcal{B}^+(y)(e_u^-(\lambda)+\theta_1^+(y,\lambda))\md y. \label{eqtheta1plus} 
\end{split}
\eqq
Using the fact that $\eta_--\alpha>0$, we deduce that $\left\|\theta_1^+(x,\lambda)\right\|=\mathcal{O}\left(e^{-\alpha x}\right)$ as $x\rightarrow +\infty$ uniformly in $\lambda$. To obtain the conclusion of the lemma and a solution defined for all $x\geq0$, we flow backward the solution of \eqref{eqZphi} from $x=x_0$ to $x=0$.

We now construct $\psi_1^+(x)$.  We remark that $\phi_1^+(x)$ is unique up to scalar multiplication.  There is no such unique choice of $\psi_1^+(x)$, but we will have to construct $\psi_1^+(x)$ adequately so that $\kappa_1^+(x,\lambda)$ converges to zero with rate $e^{-\alpha x}$ as desired.  We begin as we did above, this time   transforming  (\ref{eqPsyst}) as $Z(x)=e^{\sqrt{\lambda}x}P(x)$ such that $Z(x)$ obeys
\bqq
 Z'(x)=\left(\mathcal{A}^+(\lambda)-\sqrt{\lambda}I_4\right)Z(x)+ \mathcal{B}^+(x)Z(x), \quad x\geq1.
 \label{eqZpsi}
 \eqq
Repeating the steps as above the solution of this equation can be expressed as 
\[ Z(x)=e^{(\mathcal{A}^+(\lambda)-\sqrt{\lambda}I_4)x}Z(x_0)+\int_{x_0}^x e^{(\mathcal{A}^+(\lambda)-\sqrt{\lambda}I_4)(x-y)} \mathcal{B}^+(y)Z(y)\md y, \quad x\geq1, \]
for some $x_0\geq1$. We then select  
\[ Z(x_0)=e_u^+(\lambda)-\int_{x_0}^\infty e^{-(\mathcal{A}^+(\lambda)-\sqrt{\lambda}I_4)y}\Pi_{cu}(\lambda) \mathcal{B}^+(y)Z(y)\md y, \]
so that 
\bqq \begin{split} Z(x)=&~e_u^+(\lambda)+\int_{x_0}^x e^{(\mathcal{A}^+(\lambda)-\sqrt{\lambda}I_4)(x-y)}\Pi_{ss}(\lambda) \mathcal{B}^+(y)Z(y)\md y\\
 &- \int_{x}^\infty e^{(\mathcal{A}^+(\lambda)-\sqrt{\lambda}I_4)(x-y)}\Pi_{cu}(\lambda) \mathcal{B}^+(y)Z(y)\md y. \end{split}  \label{eqZpsimap} 
 \eqq
We have the bounds
\begin{eqnarray*} \left\|e^{(\mathcal{A}^+(\lambda)-\sqrt{\lambda}I_4)x}\Pi_{ss}(\lambda)\right\|&\leq& K e^{-\eta_-x}, \quad x>0, \\
\left\|e^{(\mathcal{A}^+(\lambda)-\sqrt{\lambda}I_4)x}\Pi_{cu}(\lambda)\right\|&\leq& Ke^{-2\sqrt{\lambda}x}.  \quad x<0 
\end{eqnarray*} 
Note  that the constant $K$ is independent of $\lambda$.  While the ranges of $\Pi_{ws}(\lambda)$ and $\Pi_{wu}(\lambda)$ converge as $\lambda\to 0$, the projection  $\Pi_{cu}(\lambda)$ remains well behaved.  Note also that the center unstable solution grows slowly as $x$ decreases.  Nonetheless, this slow growth is offset by the strong decay of the matrix $\mathcal{B}^+(x)$ and we once again obtain a solution using the contraction mapping theorem on $L^{\infty}([x_0,\infty))$ for some $x_0$ sufficiently large.  

Expanding $Z(x)=e_u^+(\lambda)+\kappa_1^+(x,\lambda)$ we have 
\bqq
\begin{split} \kappa_1^+(x,\lambda)=& \int_{x_0}^x e^{(\mathcal{A}^+(\lambda)-\sqrt{\lambda}I_4)(x-y)}\Pi_{ss}(\lambda) \mathcal{B}^+(y)(e_u^+(\lambda)+\kappa_1^+(y,\lambda))\md y  \\
&- \int_{x}^\infty e^{(\mathcal{A}^+(\lambda)-\sqrt{\lambda}I_4)(x-y)}\Pi_{cu}(\lambda) \mathcal{B}^+(y)(e_u^+(\lambda)+\kappa_1^+(y,\lambda))\md y. 
\end{split}
 \label{eqkappa1plus}
 \eqq
We then have that $\left\|\kappa_1^+(x,\lambda)\right\|=\mathcal{O}\left(e^{-\alpha x}\right)$ as $x\rightarrow +\infty$ uniformly in $\lambda$, and this concludes the proof for $\psi_1^+(x)$ after flowing backward the solution of \eqref{eqZpsi} from $x=x_0$ to $x=0$.  

The proofs for $\phi_2^+(x)$ and $\psi_2^+(x)$ follow along similar lines. 

\end{Proof} 

The matrix $\mathcal{A}^-(\lambda)$ also has four eigenvalues:
\[ \mu_u^\pm(\lambda)=-\delta -\gamma_*\pm \sqrt{\gamma_*^2+1+\lambda}, \quad \mu_v^\pm(\lambda)=-\delta -\frac{\gamma_*}{\sigma}\pm \frac{1}{\sigma}\sqrt{\gamma_*^2+\sigma r(b-1)+\sigma \lambda}, \]
and four eigenvectors
\[ \epsilon_u^\pm(\lambda) =\left(\begin{array}{c} 1 \\ \mu_u^\pm(\lambda) \\ 0 \\ 0 \end{array}\right), \quad \epsilon_v^\pm(\lambda) =\left(\begin{array}{c} x_u^\pm(\lambda) \\ \mu_v^\pm(\lambda) x_u^\pm(\lambda)   \\ 1 \\ \mu_v^\pm(\lambda) \end{array}\right), \]
with
\bqs
x_u^\pm(\lambda)=\frac{a}{(\mu_v^\pm(\lambda))^2+\mu_v^\pm(\lambda)(c_*+2\delta)-\lambda+(-1+c_*\delta+\delta^2)}.
\eqs

Similarly to Lemma~\ref{lem:+solsp}, we characterize solutions of system \eqref{eqPsyst} with prescribed exponential growth and decay rates at $-\infty$.

\begin{lem} \label{lem:-solsm} For $x\leq 0$ and for all $|\lambda|<M_s$, to the right of $\Gamma$ and off the negative real axis, we have the existence of two linearly independent bounded solutions $\phi_{1,2}^-$ together with two linearly independent unbounded solutions $\psi_{1,2}^-$ with prescribed asymptotic growth rates,
\bqq
\begin{split}
\phi_1^-(x)&= e^{\mu_u^+(\lambda) x}\left( \epsilon_u^+(\lambda)+\theta_1^-(x,\lambda)\right),  \\
 \phi_2^-(x)&= e^{\mu_v^+(\lambda) x}\left( \epsilon_v^+(\lambda)+\theta_2^-(x,\lambda)\right),  \\
\psi_1^-(x)&= e^{\mu_u^-(\lambda) x}\left( \epsilon_u^-(\lambda)+\kappa_1^-(x,\lambda)\right), \\
\psi_2^-(x)&= e^{\mu_v^-(\lambda) x}\left( \epsilon_v^-(\lambda)+\kappa_2^-(x,\lambda)\right).
\end{split} \label{eq:+sols} \eqq
Furthermore, we have the following uniform bounds for the the vectors in (\ref{eq:+sols}),
\bqs
\left\|\theta_{1,2}^-(x,\lambda)\right\|\leq C e^{\alpha x} \text{ and } \left\|\kappa_{1,2}^-(x,\lambda)\right\|\leq C e^{\alpha x}, \quad x\leq0, 
\eqs 
 for some constant $C>0$ independent of $\lambda$.
\end{lem}

\begin{Proof} This is a consequence of Proposition 3.1 of \cite{zumbrun98}.  We also refer the reader to the proof of Lemma~\ref{lem:+solsp}. 
\end{Proof} 

Note that as $\lambda\to 0$, $e_u^+(\lambda)$ and $e_v^+(\lambda)$ become colinear.  The following lemma shows that the same holds for $\phi^+_1$ and $\psi_1^+$.  It will be important in Section 3 to derive the required bounds on $G_\lambda(x,y)$.  

\begin{lem}\label{lem:diff} For $x\geq 0$ and for all $|\lambda|<M_s$, to the right of $\Gamma$ and off the negative real axis, it holds that 
\[ \theta_1^+(x,\lambda)-\kappa_1^+(x,\lambda)=\sqrt{\lambda}\Lambda(x,\lambda), \]
where $\Lambda(x,\lambda)=\O(|x|e^{-\alpha x})$ uniformly in $\lambda$, and $\Lambda$ is analytic in $\sqrt{\lambda}$. 
\end{lem}
\begin{Proof} Recall the expressions (\ref{eqtheta1plus}) and (\ref{eqkappa1plus}).  Taking their difference we find
\bqq
\begin{split} \theta_1^+(x,\lambda)-\kappa_1^+(x,\lambda) =& \int_{x_0}^x e^{\mathcal{A}^+(\lambda)(x-y)}\Pi_{ss}(\lambda) \mathcal{B}^+(y) \left(e^{\sqrt{\lambda}(x-y)}e_u^-(\lambda)- e^{-\sqrt{\lambda}(x-y)}e_u^+(\lambda)\right)\md y  \\ 
&- \int_{x}^\infty e^{\mathcal{A}^+(\lambda)(x-y)}\Pi_{cu}(\lambda) \mathcal{B}^+(y)\left(e^{\sqrt{\lambda}(x-y)}e_u^-(\lambda)- e^{-\sqrt{\lambda}(x-y)}e_u^+(\lambda)\right)\md y  \\
&+ \int_{x_0}^x e^{\mathcal{A}^+(\lambda)(x-y)}\Pi_{ss}(\lambda) \mathcal{B}^+(y) \left(e^{\sqrt{\lambda}(x-y)}\theta_1^+(y,\lambda)- e^{-\sqrt{\lambda}(x-y)}\kappa_1^+(y,\lambda)\right)\md y  \\
&- \int_{x}^\infty e^{\mathcal{A}^+(\lambda)(x-y)}\Pi_{cu}(\lambda) \mathcal{B}^+(y)\left(e^{\sqrt{\lambda}(x-y)}\theta_1^+(y,\lambda)- e^{-\sqrt{\lambda}(x-y)}\kappa_1^+(y,\lambda)\right)\md y.
\end{split}
 \label{eqthetaminuskappa}
\eqq 
It can be verified from the expression (\ref{eq:yv}) that 
\[ e^{\sqrt{\lambda}(x-y)}e_u^-(\lambda)- e^{-\sqrt{\lambda}(x-y)}e_u^+(\lambda):=\sqrt{\lambda}\mathbf{e}\left(x-y,\sqrt{\lambda}\right),\]
with $\mathbf{e}$ analytic in $\sqrt{\lambda}$. Then, we rearrange terms to get
\bqs
e^{\sqrt{\lambda}(x-y)}\theta_1^+(y,\lambda)- e^{-\sqrt{\lambda}(x-y)}\kappa_1^+(y,\lambda)=e^{\sqrt{\lambda}(x-y)}\left(\theta_1^+(y,\lambda)-\kappa_1^+(y,\lambda)\right)+2\sqrt{\lambda}\frac{\sinh(\sqrt{\lambda}(x-y))}{\sqrt{\lambda}}\kappa_1^+(y,\lambda).
\eqs
We see that we can rewrite the difference $\theta_1^+(x,\lambda)-\kappa_1^+(x,\lambda)$ as follows
\bqs
\theta_1^+(x,\lambda)-\kappa_1^+(x,\lambda)=\sqrt{\lambda} \Phi(x,\sqrt{\lambda})+\mathcal{T}\cdot \left(\theta_1^+(x,\lambda)-\kappa_1^+(x,\lambda) \right),
\eqs
where we have set
\begin{align*}
\Phi(x,\sqrt{\lambda}):=&\int_{x_0}^x e^{\mathcal{A}^+(\lambda)(x-y)}\Pi_{ss}(\lambda) \mathcal{B}^+(y)\mathbf{e}\left(x-y,\sqrt{\lambda}\right)  \md y\\
&- \int_{x}^\infty e^{\mathcal{A}^+(\lambda)(x-y)}\Pi_{cu}(\lambda) \mathcal{B}^+(y)\mathbf{e}\left(x-y,\sqrt{\lambda}\right) \md y\\
&+2\int_{x_0}^x e^{\mathcal{A}^+(\lambda)(x-y)}\Pi_{ss}(\lambda) \mathcal{B}^+(y)\frac{\sinh(\sqrt{\lambda}(x-y))}{\sqrt{\lambda}}\kappa_1^+(y,\lambda)  \md y\\
&- 2\int_{x}^\infty e^{\mathcal{A}^+(\lambda)(x-y)}\Pi_{cu}(\lambda) \mathcal{B}^+(y)\frac{\sinh(\sqrt{\lambda}(x-y))}{\sqrt{\lambda}}\kappa_1^+(y,\lambda)\md y,
\end{align*}
and the operator $\mathcal{T}_\lambda$ is defined as
\begin{align*}
\mathcal{T}_\lambda\cdot Z(x) :=& \int_{x_0}^x e^{(\mathcal{A}^+(\lambda)+\sqrt{\lambda}I_4)(x-y)}\Pi_{ss}(\lambda) \mathcal{B}^+(y)Z(y)  \md y\\
&- \int_{x}^\infty e^{(\mathcal{A}^+(\lambda)+\sqrt{\lambda}I_4)(x-y)}\Pi_{cu}(\lambda) \mathcal{B}^+(y)Z(y)\md y.
\end{align*}
Using the decaying properties of $\mathcal{B}^+$ and $\kappa_1^+$, we get that
\bqs
\left\| \Phi(x,\sqrt{\lambda}) \right\| = \mathcal{O}\left( |x|e^{-\alpha x}\right),
\eqs
as $x\rightarrow+\infty$ uniformly in $\lambda$. We also note that both $\Psi$ and the map $\mathcal{T}_\lambda$ are analytic in $\sqrt{\lambda}$. As a consequence, using an iterative argument we get the conclusion of the lemma.
\end{Proof}

\subsection{The Evans Function}

The Evans function; see \cite{alexander90},  is defined as,
\bqq \W_\lambda(0)=\mathrm{det}\left(\begin{matrix} \phi_1^+(0) &  \phi_2^+(0) &\phi_1^-(0) &\phi_2^-(0)\end{matrix} \right).\label{eq:Evans} \eqq

The following fact is essential to our analysis.
\begin{lem}\label{lem:Evans} Assume that assumptions \eqref{par_mono}-\eqref{par_lin} hold, then $\W_0(0)\neq 0$.  
\end{lem}
\begin{Proof}
The extension of the Evans function to $\lambda=0$ is possible due to the Gap Lemma; see \cite{gardner98, kapitula98}.  That $\W_0(0)\neq 0$ is a direct consequence of the fact that $\lambda=0$ is not an eigenvalue of $\cL$ and so does not contribute to a zero of the Evans function $\W_0$. This is proved in Lemma~\ref{lem:0eig} using comparison principle techniques.
\end{Proof}

Define the following subset of the complex plane,
\bqq \Omega_\delta=\left\{ \lambda\in\mathbb{C}\ \  | \  \mathrm{Re}(\lambda)\geq -\delta_0-\delta_1 |\mathrm{Im}(\lambda)| \right\} \label{eqOmegadelta}, \eqq
for $\delta_0$ and $\delta_1$ to be specified below.  We denote the boundary of this region by $\Gamma_\delta$.  

\begin{lem}\label{lemnozeros} There exists $\delta_0>0$ and $\delta_1>0$ such that $\W_\lambda(0)\neq 0$ for all $\lambda \in \Omega_\delta$ and off the negative real axis.
\end{lem}
\begin{Proof} We recall two important properties of the Evans function $\W_\lambda(0)$.  First, $\W_\lambda(0)$ is analytic in $\lambda$.  Second, there exists an $M>0$, large enough, such that for all  $|\lambda|>M$ and to the right of the essential spectrum we have that $\W_\lambda(0)\neq 0$. Next, we apply Lemma~\ref{lem:noposeig} which says that there is no positive eigenvalues of $\cL$ such that we get that $\W_\lambda(0)\neq 0$ for all $\lambda$ with $\Re(\lambda)\geq0$. Here, we have also used the fact that the operator $\cL$ is monotone due to the competitive structure of the Lotka-Volterra system which implies that any eigenvalue of $\cL$ has to be real. In particular, we can actually conclude that $\sigma_{\mathrm{pt}}(\cL)=\emptyset$ and that the essential spectrum is contained in some sector, given by a parabola, in the left-half complex plane at the exception of the negative real axis. As a consequence, it is possible to choose $\delta_0$ and $\delta_1$ so that the result holds.  
\end{Proof}

In the next section, we will require expressions for $\W_\lambda(x)$ for nonzero $x$.  We recall that, 
\bqs 
\W_\lambda(x)=\mathrm{det}\left(\begin{matrix} \phi_1^+(x) &  \phi_2^+(x) &\phi_1^-(x) &\phi_2^-(x)\end{matrix} \right), \quad \forall x\in\R,
\]
is a Wronskian of linearly independent solutions to (\ref{eqPsyst}). As a consequence, we have that 
\[ \frac{\md}{\md x} \W_\lambda(x)=-\left( c_*+2\frac{\omega(x)'}{\omega(x)}+\frac{c_*}{\sigma} +2\sigma \frac{\omega'(x)}{\omega(x)}\right)\W_\lambda(x).\] 
We therefore compute that for $x\geq 1$
\bqq \W_\lambda(x)=\W_\lambda(0) e^{(\nu_v^+(\lambda)+\nu_v^-(\lambda))x}, \label{eq:+W} \eqq
whereas for $x\leq -1$, 
\bqq \W_\lambda(x)=\W_\lambda(0) e^{(\mu_u^+(\lambda)+\mu_u^-(\lambda)+\mu_v^+(\lambda)+\mu_v^-(\lambda))x}. \label{eq:-W} \eqq

\section{The pointwise Green's function}

For $\lambda\in\C$, the pointwise Green's function $\bG_\lambda(x,y)$ is a $2\times2$ matrix satisfying 
\bqq \left(\mathcal{L}-\lambda \mathrm{I}_2\right)  \bG_\lambda(x,y)=-\delta(x-y)\mathrm{I}_2, \quad \forall x,y\in\R, \label{eqGdef} \eqq
whose components are defined as
\[ \bG_\lambda(x,y):=\left(\begin{array}{cc} \bG_\lambda^{11}(x,y) & \bG_\lambda^{12} (x,y) \\
\bG_\lambda^{21}(x,y) & \bG_\lambda^{22}(x,y) \end{array}\right). \]

The main result of this section is as follows. We recall that $M_s>0$ is defined such that for all $|\lambda|<M_s$, to the right of $\Gamma$ and off the negative real axis, the ordering \eqref{eq:orderingeig} of the eigenvalues of the asymptotic matrix $\mathcal{A}^+(\lambda)$ is valid.

\begin{prop}\label{thmGlambda} The pointwise Green's function $\bG_\lambda(x,y)$ satisfies the following bounds:
\begin{itemize}
\item[(i)]For all $\lambda$ to the right of $\Gamma$ off the imaginary axis and with $|\lambda|<M_s$ it holds that  
\bqq 
\bG_\lambda(x,y)= \mathbf{H}_\lambda(x,y) e^{-\sqrt{\lambda}|x-y|}, \quad \forall x,y\in\R,
\label{eqGlambdanearzero}
\eqq
where $\mathbf{H}_\lambda$ is an analytic function of $\lambda$, bounded uniformly as a function of $x$ and $y$.  
\item[(ii)] There exist $M_l>0$, $C_l>0$ and $\eta>0$ such that for all $\lambda$ to the right of $\Gamma$ and with $|\lambda|>M_l$ it holds that 
\bqq 
\left| \bG_\lambda^{ij}(x,y)\right| \leq \frac{C_l}{\sqrt{|\lambda|}} e^{-\sqrt{|\lambda|}\eta  |x-y|}, \quad \forall x,y\in\R,
\label{eqGlambdabig}
\eqq
for all $i,j\in\left\{1,2\right\}$.
\item[(iii)] For all $\lambda$ to the right of $\Gamma$ and with $M_s\leq |\lambda|\leq M_l$ there exists a $C_m>0$ such that 
\bqq \left| \bG_\lambda^{ij}(x,y)\right| \leq C_m, \quad \forall x,y\in\R, \label{eqGlambdamedium} \eqq
for all $i,j\in\left\{1,2\right\}$.
\end{itemize}
\end{prop}

The proof of the large $\lambda$ estimate \eqref{eqGlambdabig} relies on a rescaling of the independent variable so that (\ref{eqGdef}) becomes independent of $x$ to leading order in $\frac{1}{|\lambda|}$.  Since this is a standard result we refer the reader to  Proposition 7.3 of \cite{zumbrun98}. For $\lambda$ to the right of $\Gamma$ and with $M_s\leq |\lambda|\leq M_l$, there remains a spectral gap between the stable eigenvalues of the asymptotic system at $+\infty$ and the unstable eigenvalues of the asymptotic system at $-\infty$.  When combined with the fact that there are no eigenvalues in this region; see Lemma~\ref{lemnozeros}, we find that the first order system has a generalized exponential dichotomy on the whole real line and boundedness of the second order Green's function $\bG_\lambda(x,y)$ follows. 

In the remaining of this section, we establish the remaining estimate \eqref{eqGlambdanearzero} of Proposition \ref{thmGlambda}. The general idea is to derive expressions for the pointwise Green's function $\bG_\lambda$ in terms of the bounded solutions $\phi_{1,2}^\pm$ defined in Lemma~\ref{lem:+solsp} and Lemma~\ref{lem:-solsm}.

\subsection{General formulation} 

Notation: Define the projections $\pi_j:\mathbb{C}^4\to \mathbb{C}^3$, whereby the $j$th component is deleted.  Similarly, define $\pi_{ij}:\mathbb{C}^4\to \mathbb{C}^2$ where the $i$th and $j$th component are deleted.  In an analogous manner define $\pi_{(i)^c}:\mathbb{C}^4\to \mathbb{C}^1$ where all but the $i$th component is removed. 

Throughout, we assume that $\lambda$ is to the right of $\Gamma$ off the imaginary and with $|\lambda|<M_s$ such that the estimates derived in Lemma~\ref{lem:+solsp} and Lemma~\ref{lem:-solsm} are satisfied.

For $x\neq y$, the relation (\ref{eqGdef}) indicates that the components should obey
\begin{align*}
(\mathcal{L}_u-\lambda) \bG_\lambda^{11}(x,y)+\mathcal{L}_{12} \bG_\lambda^{21}(x,y)&= 0, \\
(\mathcal{L}_u-\lambda) \bG_\lambda^{12}(x,y)+\mathcal{L}_{12} \bG_\lambda^{22}(x,y)&= 0, \\
\mathcal{L}_{21} \bG_\lambda^{11}(x,y)+(\mathcal{L}_v-\lambda) \bG_\lambda^{21}(x,y)&= 0, \\
\mathcal{L}_{21} \bG_\lambda^{12}(x,y)+(\mathcal{L}_v-\lambda) \bG_\lambda^{22}(x,y)&= 0.
\end{align*}
Thus, we require bounded solutions for the system of equations
\begin{align*}
(\mathcal{L}_u-\lambda) \bH+\mathcal{L}_{12} \bK&= 0 \\
\mathcal{L}_{21} \bH+(\mathcal{L}_v-\lambda) \bK&= 0,
\end{align*}
from which we can use $(\bH,\bK)=(\bG_\lambda^{11},\bG_\lambda^{21}) $ or $(\bH,\bK)=(\bG_\lambda^{12},\bG_\lambda^{22}) $ to construct the Green's function.  Decay requires that for $x>y$
\[ \left(\begin{array}{c} \bH \\ \bK \end{array}\right)(x,y) =c_1(y)\pi_{2,4} \phi_1^+(x) +c_2(y) \pi_{2,4} \phi_2^+(x), \]
for some $c_1(y)$ and $c_2(y)$ while for $x<y$ we require 
\[ \left(\begin{array}{c} \bH \\ \bK \end{array}\right)(x,y) =-c_3(y)\pi_{2,4} \phi_1^-(x) -c_4(y) \pi_{2,4} \phi_2^-(x), \]
for some $c_3(y)$ and $c_4(y)$. We now impose matching conditions at $x=y$.  For the pair  $(\bG_\lambda^{11},\bG_\lambda^{21}) $ we require $\bG_\lambda^{21}$ to be continuously differentiable while we require $\bG_\lambda^{11}$ to be continuous with a jump discontinuity in the derivative.  Thus to determine the coefficients $c_j(y)$ in this case leads to a solvability condition
\[ \left(\begin{array}{cccc} \phi_1^+(y) &  \phi_2^+(y) &\phi_1^-(y) &\phi_2^-(y)\end{array} \right) \left(\begin{array}{c} c_1(y) \\ c_2(y) \\ c_3(y) \\ c_4(y) \end{array}\right)= \left(\begin{array}{c} 0 \\ -1 \\ 0 \\ 0 \end{array}\right).\]
On the other hand, to solve for the pair $(\bG_\lambda^{12},\bG_\lambda^{22}) $ we require a jump discontinuity for $\bG_\lambda^{22}$ and we obtain conditions 
\[ \left(\begin{array}{cccc} \phi_1^+(y) &  \phi_2^+(y) &\phi_1^-(y) &\phi_2^-(y)\end{array} \right) \left(\begin{array}{c} c_1(y) \\ c_2(y) \\ c_3(y) \\ c_4(y) \end{array}\right)= \left(\begin{array}{c} 0 \\ 0 \\ 0 \\ -\frac{1}{\sigma} \end{array}\right).\]

\subsection{ The components $\bG_\lambda^{11}(x,y)$ and $\bG_\lambda^{21}(x,y)$ }\label{sec:G11}
To compute $\bG_\lambda^{11}(x,y)$ and $\bG_\lambda^{21}(x,y)$ we solve
\[  \left(\begin{array}{c} c_1(y) \\ c_2(y) \\ c_3(y) \\ c_4(y) \end{array}\right)= \left(\begin{array}{cccc} \phi_1^+(y) &  \phi_2^+(y) &\phi_1^-(y) &\phi_2^-(y)\end{array} \right)^{-1} \left(\begin{array}{c} 0 \\ -1 \\ 0 \\ 0 \end{array}\right),\]
from which we obtain
\begin{align*} 
c_1(y)&= \frac{1}{\W_\lambda(y)}\mathrm{det}\left(\begin{array}{ccc} \pi_2 \phi_2^+(y) & \pi_2 \phi_1^-(y) & \pi_2\phi_2^-(y)\end{array}\right), \\
c_2(y)&= -\frac{1}{\W_\lambda(y)}\mathrm{det}\left(\begin{array}{ccc} \pi_2 \phi_1^+(y) & \pi_2 \phi_1^-(y) & \pi_2\phi_2^-(y)\end{array}\right),\\
c_3(y)&= \frac{1}{\W_\lambda(y)}\mathrm{det}\left(\begin{array}{ccc} \pi_2 \phi_1^+(y) & \pi_2 \phi_2^+(y) & \pi_2\phi_2^-(y)\end{array}\right),\\
c_4(y)&= -\frac{1}{\W_\lambda(y)}\mathrm{det}\left(\begin{array}{ccc} \pi_2 \phi_1^+(y) & \pi_2 \phi_2^+(y) & \pi_2\phi_1^-(y)\end{array}\right),
\end{align*}
where we recall the definition of $\W_\lambda(y)$ in (\ref{eq:Evans}).  Upon introducing the following notation
\[ \phi_{j}^\pm(x)=\left(\begin{array}{c} \phi_{j,u}^\pm(x) \\ \tilde{\phi}_{j,u}^\pm(x) \\ \phi_{j,v}^\pm(x) \\ \tilde{\phi}_{j,v}^\pm(x) \end{array}\right), \quad j=1,2,  \]
for the components of $\phi_{1,2}^\pm$, we then see that 
\begin{eqnarray} \bG_\lambda^{11}(x,y)&=&\left\{ \begin{array}{cc} \frac{N_1^+(y)\phi_{1,u}^+(x)}{\W_\lambda(y) }+ \frac{N_2^+(y)\phi_{2,u}^+(x)}{\W_\lambda(y) },   & x>y, \\
\frac{N_1^-(y)\phi_{1,u}^-(x)}{\W_\lambda(y) }+ \frac{N_2^-(y)\phi_{2,u}^-(x)}{\W_\lambda(y) },   & x<y, \end{array}\right. \label{eq:G11} \\
 \bG_\lambda^{21}(x,y)&=& \left\{ \begin{array}{cc} \frac{N_1^+(y)\phi_{1,v}^+(x)}{\W_\lambda(y) }+ \frac{N_2^+(y)\phi_{2,v}^+(x)}{\W_\lambda(y) },   & x>y, \\
\frac{N_1^-(y)\phi_{1,v}^-(x)}{\W_\lambda(y) }+ \frac{N_2^-(y)\phi_{2,v}^-(x)}{\W_\lambda(y) },   & x<y .\end{array}\right. \label{eq:G21} \end{eqnarray}
In both of these expressions,
\begin{eqnarray*} N_1^+(y)&=&\mathrm{det}\left(\begin{array}{ccc} \pi_2 \phi_2^+(y) & \pi_2 \phi_1^-(y) & \pi_2\phi_2^-(y)\end{array}\right), \\
N_2^+(y)&=&- \mathrm{det}\left(\begin{array}{ccc} \pi_2 \phi_1^+(y) & \pi_2 \phi_1^-(y) & \pi_2\phi_2^-(y)\end{array}\right),\\
N_1^-(y)&=& -\mathrm{det}\left(\begin{array}{ccc} \pi_2 \phi_1^+(y) & \pi_2 \phi_2^+(y) & \pi_2\phi_2^-(y)\end{array}\right), \\
N_2^-(y)&=&  \mathrm{det}\left(\begin{array}{ccc} \pi_2 \phi_1^+(y) & \pi_2 \phi_2^+(y) & \pi_2\phi_1^-(y)\end{array}\right) .
\end{eqnarray*}

We introduce some further notation before proceeding.  Define $\bm$ and $\um$ with 
\[ \mathrm{Re}(\um) \leq \mathrm{Re}(\bm)  \leq 0 \]
as
\begin{eqnarray*}
\bm &:=& \left\{ \begin{array}{ccc} \mu_u^-(\lambda)  &\text{ if }& \mathrm{Re}\left( \mu_u^-(\lambda)-\mu_v^-(\lambda)\right)>0, \\
\mu_v^-(\lambda)  & \text{ if }& \mathrm{Re}\left( \mu_u^-(\lambda)-\mu_v^-(\lambda)\right)\leq 0,
\end{array}\right.  \\
\um &:=& \left\{ \begin{array}{ccc} \mu_v^-(\lambda)  &\text{ if } & \mathrm{Re}\left( \mu_u^-(\lambda)-\mu_v^-(\lambda)\right)>0, \\
\mu_u^-(\lambda)  &\text{ if } & \mathrm{Re}\left( \mu_u^-(\lambda)-\mu_v^-(\lambda)\right)\leq 0.
\end{array}\right. 
\end{eqnarray*}

\begin{lem}\label{lem:Nbds}For all $\lambda$ to the right of $\Gamma$ off the imaginary axis and with $|\lambda|<M_s$ it holds:
\begin{itemize}
\item[(i)] For $y>0$:
\begin{eqnarray*} 
N_1^+(y)&=&  e^{\sqrt{\lambda}y} e^{(\nu_v^+(\lambda)+\nu_v^-(\lambda))y} \O(1),\nonumber \\
N_2^+(y)&=&   e^{\nu_v^+(\lambda)y}\O(1), \nonumber \\
N_1^-(y)&=&  e^{-\sqrt{\lambda}y} e^{(\nu_v^+(\lambda)+\nu_v^-(\lambda))y}\O(1), \nonumber \\
N_2^-(y)&=&   e^{-\sqrt{\lambda}y} e^{(\nu_v^+(\lambda)+\nu_v^-(\lambda))y}\O(1).
\end{eqnarray*}
\item[(ii)] For $y<0$:
\begin{eqnarray*} 
N_1^+(y)&=&  e^{(\mu_u^+(\lambda)+\mu_v^+(\lambda)+\um)y} \O(1),\nonumber \\
N_2^+(y)&=&  e^{(\mu_u^+(\lambda)+\mu_v^+(\lambda)+\um )y}\O(1), \nonumber \\
N_1^-(y)&= &  e^{(\mu_v^+(\lambda)+\mu_v^-(\lambda)+\mu_u^-(\lambda)) )y}\O(1), \nonumber \\
  N_2^-(y)&=&  e^{(\mu_u^+(\lambda)+\mu_v^-(\lambda)+\mu_u^-(\lambda)) )y} \O(1) .
\end{eqnarray*}
\end{itemize}
In both cases, $\O(1)$ refer to terms that are analytic in $\lambda$ and bounded uniformly in $y$.
\end{lem}

\begin{Proof} We focus on the first case: $N_1^+(y)$ for $y>0$.  Note that an explicit expression for $\phi_2^+(y)$ is provided in Lemma~\ref{lem:+solsp}.  Since $y>0$, we do not have a similar representations  for $\phi_{1,2}^-(y)$.  We now go about computing bounds.  We need estimates for 
\bqq \mathrm{det}\left(\begin{array}{cc}  \pi_{2,j} \phi_1^-(y) & \pi_{2,j} \phi_2^-(y)\end{array}\right),\label{eq:detex} \eqq
for $j={1,3,4}$.  First note that $\phi_1^-(y)$ and $\phi_2^-(y)$ are linearly independent by construction and so the above determinant is non zero for at least some choice of $j$.  To understand how this determinant evolves, we interpret (\ref{eq:detex}) in the language of differential forms; see for example \cite{jones94}.  Let $\eta_i$ denote the differential 1-form, $\eta_i=\pi_{(i)^c}$, which when applied to a vector extracts the $i$th component.  The differential 2-form $\eta_{ij}=\eta_i \wedge \eta_j$ acts on a pair of vectors,
\[ \eta_{ij}(\phi,\psi)=\mathrm{det}\left(\begin{array}{cc} \eta_i(\phi) & \eta_i (\psi) \\
\eta_j(\phi) & \eta_j(\psi) \end{array}\right).\]
Thus, the expression in (\ref{eq:detex}) can also be written as $\eta_{(2j)^c}(\phi_1^-(y),\phi_2^-(y))$.  We want to understand how this quantity evolves in $y$.  

Note that $\phi_{1,2}^-(y)$ are solutions of (\ref{eqPsyst}) and  the one form $\eta_i(\phi_{1,2}^-)$ is just the $i$th component of the solution.  Due to the exponential convergence of $\mathcal{B}(y)$ we obtain the bounds $|\eta_i(\phi_{j}^-(y))|\leq C e^{\nu_v^+(\lambda)y}$.  For the two forms, it is convenient to first diagonalize (\ref{eqPsyst}) via the transformation $Q=SP$, from which 
\bqq
Q'=\Theta(y,\lambda)Q+\mathcal{C}(y)Q, 
\label{eqQsyst}
\eqq
where $\Theta=\mathrm{diag}(\nu_v^+(\lambda),\sqrt{\lambda},-\sqrt{\lambda},\nu_v^-(\lambda))$ and $\mathcal{C}(y)=S\mathcal{B}(y) S^{-1}$.  Once diagonalized, it is straightforward to derive differential equations for the two forms via the identity
\[ \frac{\md}{\md y} \eta_{ij}(\phi(y),\psi(y)) = \frac{\md}{\md y}\left( \eta_i(\phi(y))\eta_j (\psi(y)) -\eta_i(\psi(y))\eta_j(\phi(y))\right), \]
and we obtain a six dimensional system of ODEs for $\Xi(y)$ describing the evolution of the six non-trivial two forms,
\bqq \Xi'=\Upsilon \Xi+\mathcal{D}(y)\Xi,\label{eq:Xi} \eqq
where $\Upsilon=\mathrm{diag}(\nu_v^+(\lambda)+\sqrt{\lambda},\nu_v^+(\lambda)-\sqrt{\lambda},\nu_v^+(\lambda)+\nu_v^-(\lambda),0, \sqrt{\lambda}-\nu_v^-(\lambda),\sqrt{\lambda}-\nu_v^-(\lambda))$. Once again, $\mathcal{D}(y)$ converges exponentially and therefore any solution of (\ref{eq:Xi}) satisfies $|\Xi(y)|\leq C e^{(\nu_v^+(\lambda)+\sqrt{\lambda})y}$.  Finally, since (\ref{eq:detex}) is constructed from linear combinations of the components of $\Xi$ we obtain the same bound there.  

\end{Proof}

We are now able to use the formulas (\ref{eq:G11}) and (\ref{eq:G21}) together with the bounds from Lemma~\ref{lem:Nbds} to obtain the following result.  
\begin{lem}\label{lem:G11bds} Assume that $\lambda$ is to the right of $\Gamma$ off the imaginary axis and with $|\lambda|<M_s$, then we have the following estimates.
\begin{itemize}
\item[(i)] For $y\leq 0 \leq x$:
\begin{eqnarray*}
\bG_\lambda^{11}(x,y)&= &\O(1) e^{-\sqrt{\lambda}x -\bm y }+\O(1) e^{\nu_v^-(\lambda)x -\bm y } ,  \\
 \bG_\lambda^{21}(x,y)&= &\O(1) e^{-\sqrt{\lambda}x -\bm y } +  \O(1) e^{\nu_v^-(\lambda)x -\bm y }.
\end{eqnarray*}
\item[(ii)] For $x\leq 0 \leq y$:
\begin{eqnarray*}
\bG_\lambda^{11}(x,y)&= & \O(1)e^{\mu_u^+(\lambda)x -\sqrt{\lambda}y}+\O(1)e^{\mu_v^+(\lambda)x -\sqrt{\lambda}y},\\
 \bG_\lambda^{21}(x,y)&= & \O(1)e^{\mu_u^+(\lambda)x -\sqrt{\lambda}y}+\O(1)e^{\mu_v^+(\lambda)x -\sqrt{\lambda}y}.
\end{eqnarray*}
\item[(iii)]For $0\leq y \leq x$:
\begin{eqnarray*}
\bG_\lambda^{11}(x,y)&= & \O(1)e^{-\sqrt{\lambda}(x-y)}+\O(1)e^{\nu_v^-(\lambda)(x-y)-\alpha x}, \\
 \bG_\lambda^{21}(x,y)&= & \O(1)e^{-\sqrt{\lambda}(x-y)-\alpha x}+\O(1)e^{\nu_v^-(\lambda)(x-y)}.
\end{eqnarray*}
\item[(iv)]For $0\leq x \leq y$:
\begin{eqnarray*}
\bG_\lambda^{11}(x,y)&= & \O(1)e^{-\sqrt{\lambda}(y-x)}, \\
 \bG_\lambda^{21}(x,y)&= & \O(1)e^{-\sqrt{\lambda}(y-x)}.
\end{eqnarray*}
\item[(v)]For $y\leq x \leq 0$:
\begin{eqnarray*}
\bG_\lambda^{11}(x,y)&= & \O(1) e^{\mu_v^-(\lambda)(x-y)}+\O(1)e^{\mu_u^-(\lambda)(x-y)}, \\
 \bG_\lambda^{21}(x,y)&= &  \O(1) e^{\mu_v^-(\lambda)(x-y)}+\O(1)e^{\mu_u^-(\lambda)(x-y)}.
\end{eqnarray*}
\item[(vi)]For $x\leq y \leq 0$:
\begin{eqnarray*}
\bG_\lambda^{11}(x,y)&= & \O(1)e^{\mu_u^+(\lambda)(x-y)}+\O(1)e^{\mu_v^+(\lambda)(x-y)}, \\
 \bG_\lambda^{21}(x,y)&= & \O(1)e^{\mu_u^+(\lambda)(x-y)}+\O(1)e^{\mu_v^+(\lambda)(x-y)},
\end{eqnarray*}
\end{itemize}
where $\O(1)$ refer to terms that are analytic in $\lambda$ and bounded uniformly in $(x,y)$.  
\end{lem}
\begin{Proof} The result follows by combining Lemma~\ref{lem:+solsp}, Lemma~\ref{lem:-solsm} and Lemma~\ref{lem:Nbds} with the formulas (\ref{eq:+W}) and (\ref{eq:-W}). Some cases are more straightforward than others.  We comment on several now.

Cases (i) and (iii): Since $x>0$ we have expressions for $\phi_{1}^+(x)$ and $\phi_{2}^+(x)$ directly from Lemma~\ref{lem:+solsp}.  Combining these estimates with those from Lemma~\ref{lem:Nbds} we obtain the desired estimates.

Case (ii): The estimates follow directly from a direct application of Lemmas \ref{lem:+solsp}, \ref{lem:-solsm} and \ref{lem:Nbds}.

Case (iv): To obtain the bound for $\bG^{11}_\lambda(x,y)$ in this example, we recall the expression for $x<y$,
\[ \bG^{11}_\lambda(x,y)=\frac{N_1^-(y)\phi_{1,u}^-(x)}{\W_\lambda(y) }+ \frac{N_2^-(y)\phi_{2,u}^-(x)}{\W_\lambda(y) }.\]
We then expand $\cW_\lambda^{11}(x,y):=N_1^-(y)\phi_{1,u}^-(x) +N_2^-(y)\phi_{2,u}^-(x)$ into
\bqq 
\begin{split}
\cW_\lambda^{11}(x,y)=&-\mathrm{det}\left(\begin{matrix} \pi_2 \phi_1^+(y) & \pi_2 \phi_2^+(y) & \pi_2\phi_2^-(y)\end{matrix}\right)\phi_{1,u}^-(x)\\
 &+ \mathrm{det}\left(\begin{matrix} \pi_2 \phi_1^+(y) & \pi_2 \phi_2^+(y) & \pi_2\phi_1^-(y)\end{matrix}\right) \phi_{2,u}^-(x).
 \end{split}
 \label{eqdetforG11caseiv} \eqq 
Since $x>0$, we do not have explicit bounds on  $\phi_1^-(x)$ and $\phi_2^-(x)$. Instead we express these solutions as
\bqq
\begin{split}
\phi_1^-(x)&= A_1^+\phi_1^+(x)+B_1^+\psi_1^+(x)+A_2^+\phi_2^+(x) +B_2^+\psi_2^+(x),  \\
\phi_2^-(x)&= C_1^+\phi_1^+(x)+D_1^+\psi_1^+(x)+C_2^+\phi_2^+(x) +D_2^+\psi_2^+(x).
\end{split}
 \label{eqAsBsCsDs}
\eqq
The dependence of $A_{1,2}^+$, $B_{1,2}^+$, $C_{1,2}^+$ and $D_{1,2}^+$ on $\lambda$ is suppressed but important.  Expressions for these quantities are available using Cramer's Rule.  For example, we find that 
\[ B_2^+(\lambda)=\frac{\mathrm{det}\left(\begin{matrix}  \phi_1^+(x) & \psi_1^+(x) & \phi_2^+(x) &  \phi_1^-(x) \end{matrix}\right)}{\mathrm{det}\left(\begin{matrix}  \phi_1^+(x) &  \psi_1^+(x) & \phi_2^+(x) &  \psi_2^+(x) \end{matrix}\right)},\]
while 
\[ B_1^+(\lambda)=\frac{\mathrm{det}\left(\begin{matrix}  \phi_1^+(x) & \phi_1^-(x) & \phi_2^+(x) &  \psi_2^+(x) \end{matrix}\right)}{\mathrm{det}\left(\begin{matrix}  \phi_1^+(x) &  \psi_1^+(x) & \phi_2^+(x) &  \psi_2^+(x) \end{matrix}\right)}.\]
Since both determinants are Wronskians of a set of linearly independent solutions the $x$ dependence of each cancels and $B_2^+(\lambda)$ is independent of $x$.   

First, note that $A_2^+$, $B_2^+$, $C_2^+$ and $D_2^+$ are all $\O(1)$ in $\lambda$ while $A_1^+$, $B_1^+$, $C_1^+$ and $D_1^+$ are $\O(\lambda^{-1/2})$.  We now return to (\ref{eqdetforG11caseiv}). Substituting the expansions (\ref{eqAsBsCsDs}) into the determinants in (\ref{eqdetforG11caseiv}) the expression reduces to 
\begin{align*}
\mathrm{det}\left(\begin{matrix} \pi_2 \phi_1^+(y) & \pi_2 \phi_2^+(y) & \pi_2\phi_2^-(y)\end{matrix}\right)=&~ D_1^+\mathrm{det}\left(\begin{matrix} \pi_2 \phi_1^+(y) & \pi_2 \phi_2^+(y) & \pi_2\psi_1^+(y)\end{matrix}\right)\\
&+D_2^+\mathrm{det}\left(\begin{matrix} \pi_2 \phi_1^+(y) & \pi_2 \phi_2^+(y) & \pi_2\psi_2^+(y)\end{matrix}\right),\\
\mathrm{det}\left(\begin{matrix} \pi_2 \phi_1^+(y) & \pi_2 \phi_2^+(y) & \pi_2\phi_1^-(y)\end{matrix}\right)=&~ B_1^+\mathrm{det}\left(\begin{matrix} \pi_2 \phi_1^+(y) & \pi_2 \phi_2^+(y) & \pi_2\psi_1^+(y)\end{matrix}\right)\\
&+B_2^+\mathrm{det}\left(\begin{matrix} \pi_2 \phi_1^+(y) & \pi_2 \phi_2^+(y) & \pi_2\psi_2^+(y)\end{matrix}\right).
\end{align*}
As a consequence, upon denoting
\bqs
\D_j(y):=\mathrm{det}\left(\begin{matrix} \pi_2 \phi_1^+(y) & \pi_2 \phi_2^+(y) & \pi_2\psi_j^+(y)\end{matrix}\right), \quad j=1,2,
\eqs
we get that the expression \eqref{eqdetforG11caseiv} rewrites as
\begin{align*}
\cW_\lambda^{11}(x,y)=&~\D_1(y)\left( \phi_{1,u}^+(x)(B_1^+C_1^+-D_1^+A_1^+) +\phi_{2,u}^+(x)(B_1^+C_2^+-D_1^+A_2^+)+\psi_{2,u}^+(x)(B_1^+D_2^+-D_1^+B_2^+)\right)\\
&+\D_2(y)\left( \phi_{1,u}^+(x)(B_2^+C_1^+-D_2^+A_1^+) +\psi_{1,u}^+(x)(B_2^+D_1^+-D_2^+B_1^+)+\phi_{2,u}^+(x)(B_2^+C_2^+-D_2^+A_2^+)\right) .
\end{align*}

We can now begin to estimate the first group of terms that appears in $\bG^{11}_\lambda(x,y)$, that is
\bqs 
\frac{\D_1(y)}{\W_\lambda(y) } \left( \phi_{1,u}^+(x)(B_1^+C_1^+-D_1^+A_1^+) +\phi_{2,u}^+(x)(B_1^+C_2^+-D_1^+A_2^+)+\psi_{2,u}^+(x)(B_1^+D_2^+-D_1^+B_2^+)\right) 
\eqs
and obtain that this expression is dominated by the last term in the above expression and given as
\[ e^{-\nu_v^+(\lambda)(y-x)}\O(1). \]
Here we have used the fact that for $y>0$,
\bqs
\D_1(y)=\sqrt{\lambda}e^{\nu_v^-(\lambda)y}\O(1), \quad \W_\lambda(y)= e^{\left(\nu_v^+(\lambda)+\nu_v^-(\lambda)\right)y}\O(1),
\eqs
and
\bqs
\psi_{2,u}^+(x)(B_1^+D_2^+-D_1^+B_2^+)=\frac{1}{\sqrt{\lambda}}e^{\nu_v^+(\lambda)x}\O(1).
\eqs

In a similar fashion, a term by term analysis of the second group of terms in $\bG^{11}_\lambda(x,y)$
\bqs
\frac{\D_2(y)}{\W_\lambda(y)}\left( \phi_{1,u}^+(x)(B_2^+C_1^+-D_2^+A_1^+) +\psi_{1,u}^+(x)(B_2^+D_1^+-D_2^+B_1^+)+\phi_{2,u}^+(x)(B_2^+C_2^+-D_2^+A_2^+)\right)
\eqs
yields a naive bound of
\[ \frac{1}{\sqrt{\lambda}}e^{-\sqrt{\lambda}(y-x)}\O(1). \]
Unfortunately, this bound will be insufficient to obtain the algebraic decay rates that we desire, and in the following we show how one can improve this bound and obtain the desired estimate $e^{-\sqrt{\lambda}(y-x)}\O(1)$. We will consider bounds on the sum
\[ \phi_{1,u}^+(x)(B_2^+C_1^+-D_2^+A_1^+) +\psi_{1,u}^+(x)(B_2^+D_1^+-D_2^+B_1^+). \]
We begin with 
\[ \phi_{1,u}^+(x)C_1^++\psi_{1,u}^+(x)D_1^+,\]
and recall
\[ C_1^+(\lambda)=\frac{\mathrm{det}\left(\begin{matrix}  \phi_2^-(x) & \psi_1^+(x) & \phi_2^+(x) &  \psi_2^+(x) \end{matrix}\right)}{\mathrm{det}\left(\begin{matrix}  \phi_1^+(x) &  \psi_1^+(x) & \phi_2^+(x) &  \psi_2^+(x) \end{matrix}\right)}, \ D_1^+(\lambda)=\frac{\mathrm{det}\left(\begin{matrix}  \phi_1^+(x) & \phi_2^-(x) & \phi_2^+(x) &  \psi_2^+(x) \end{matrix}\right)}{\mathrm{det}\left(\begin{matrix}  \phi_1^+(x) &  \psi_1^+(x) & \phi_2^+(x) &  \psi_2^+(x) \end{matrix}\right)}.\]
Define 
\[ \J_\lambda(x):=\mathrm{det}\left(\begin{array}{cccc}  \phi_1^+(x) &  \psi_1^+(x) & \phi_2^+(x) &  \psi_2^+(x) \end{array}\right),\]
and note that $\J_\lambda(x)=e^{\nu_v^+(\lambda)x+\nu_v^-(\lambda)x}\O(\sqrt{\lambda})$.  
Then
\begin{align*}  \phi_{1,u}^+(x)C_1^+(\lambda)&= (1+\theta_{1,u}^+(x,\lambda))\frac{\mathrm{det}\left(\begin{matrix}  \phi_2^-(x) & e_u^++\kappa_1^+(x,\lambda) & e_v^-+\theta_2^+(x,\lambda) &  e_v^++\kappa_2^+(x,\lambda) \end{matrix}\right)}{\J_\lambda(x)e^{-(\nu_v^+(\lambda)+\nu_v^-(\lambda))x}} ,\\
\psi_{1,u}^+(x)D_1^+(\lambda)&=-(1+\kappa_{1,u}^+(x,\lambda))\frac{\mathrm{det}\left(\begin{matrix}  \phi_2^-(x) & e_u^-+\theta_1^+(x,\lambda) & e_v^-+\theta_2^+(x,\lambda) &  e_v^++\kappa_2^+(x,\lambda) \end{matrix}\right)}{\J_\lambda(x)e^{-(\nu_v^+(\lambda)+\nu_v^-(\lambda))x}}.
\end{align*}
Here $\theta_{1,u}^+(x,\lambda)$ and $\kappa_{1,u}^+(x,\lambda)$ denote the first component of $\theta_{1}^+(x,\lambda)$ and $\kappa_{1}^+(x,\lambda)$ respectively. From the explicit form of $e_u^\pm(\lambda)$ and Lemma~\ref{lem:diff} which gives that $\theta_{1}^+(x,\lambda)-\kappa_{1}^+(x,\lambda)=\sqrt{\lambda}\Lambda(x,\lambda)$ with $\left\|\Lambda(x,\lambda)\right\|=\O(|x|e^{-\alpha x})$ uniformly in $\lambda$ for $x\geq0$, we obtain that
\[ (1+\theta_{1,u}^+(x,\lambda))(e_u^++\kappa_1^+(x,\lambda)) -(1+\kappa_{1,u}^+(x,\lambda))(e_u^-+\theta_1^+(x,\lambda))=\O(\sqrt{\lambda}). \]
As a consequence, we get that
\bqs
\frac{\D_2(y)}{\W_\lambda(y)}B_2^+(\phi_{1,u}^+(x)C_1^++\psi_{1,u}^+(x)D_1^+)=e^{-\sqrt{\lambda}(y-x)}\O(1).
\eqs

The same line of argument applies to $\phi_{1,u}^+(x)A_1^++\psi_{1,u}^+(x)B_1^+$, and we obtain
\bqs
\frac{\D_2(y)}{\W_\lambda(y)}D_2^+(\phi_{1,u}^+(x)A_1^++\psi_{1,u}^+(x)B_1^+)=e^{-\sqrt{\lambda}(y-x)}\O(1).
\eqs
Lastly, we observe that desired bounds on $\phi_{2,u}^+(x)(B_2^+C_2^+-D_2^+A_2^+)$ are easily obtained due to our previous remark that $A_2^+$, $B_2^+$, $C_2^+$ and $D_2^+$ are all $\O(1)$ in $\lambda$. This completes the proof of case (iv) for $\bG_{\lambda}^{11}(x,y)$, and a similar line of analysis yields identical bounds for $\bG_{\lambda}^{21}(x,y)$.  

Case (v): Here $y\leq x\leq 0$ so we will need to expand $\phi_{1}^+(x)$ and $\phi_2^+(x)$ in terms of the basis for $x<0$.  We therefore write
\bqq
\begin{split}
\phi_1^+(x)&= A_1^-\phi_1^-(x)+B_1^-\psi_1^-(x)+A_2^-\phi_2^-(x) +B_2^-\psi_2^-(x) ,  \\
\phi_2^+(x)&= C_1^-\phi_1^-(x)+D_1^-\psi_1^-(x)+C_2^-\phi_2^-(x) +D_2^-\psi_2^-(x).
\end{split} \label{eqAsBsCsDsnegative}
\eqq
Recall the expression for $y<x$,
\[ \bG^{11}_\lambda(x,y)=\frac{N_1^+(y)\phi_{1,u}^+(x)}{\W_\lambda(y) }+ \frac{N_2^+(y)\phi_{2,u}^+(x)}{\W_\lambda(y) }   .\]
Expand $\cV_\lambda^{11}(x,y):=N_1^+(y)\phi_{1,u}^+(x) +N_2^+(y)\phi_{2,u}^+(x)$ into
\bqq 
\begin{split}
\cV_\lambda^{11}(x,y)=&~\mathrm{det}\left(\begin{matrix} \pi_2 \phi_2^+(y) & \pi_2 \phi_1^-(y) & \pi_2\phi_2^-(y)\end{matrix}\right)\phi_{1,u}^+(x)\\
& -\mathrm{det}\left(\begin{matrix} \pi_2 \phi_1^+(y) & \pi_2 \phi_1^-(y) & \pi_2\phi_2^-(y)\end{matrix}\right) \phi_{2,u}^-(x).
\end{split}
\label{eqdetforG11casev} \eqq 
Using (\ref{eqAsBsCsDsnegative}) this can be expressed as 
\begin{align*}
\cV_\lambda^{11}(x,y)=&~ \mathbb{E}_1(y) \left(\phi_{1,u}^-(x)(A_1^-D_1^--B_1^-C_1^-) + \phi_{2,u}^-(x)(A_2^-D_1^--B_1^-C_2^-) 
+  \psi_{2,u}^-(x)(B_2^-D_1^--B_1^-D_2^-)\right)  \\
&+\mathbb{E}_2(y) \left(\phi_{1,u}^-(x)(A_1^-D_2^--B_2^-C_1^-) + \phi_{2,u}^-(x)(A_2^-D_2^--B_2^-C_2^-)  +\psi_{1,u}^-(x)(D_2^-B_1^--B_2^-D_1^-)\right),
\end{align*}
with the notation
\bqs
\mathbb{E}_j(y):=\mathrm{det}\left(\begin{matrix} \pi_2 \psi_j^-(y) & \pi_2 \phi_1^-(y) & \pi_2\phi_2^-(y)\end{matrix}\right), \quad j=1,2.
\eqs
The stated estimate now follows by direction calculation and Lemma~\ref{lem:-solsm}.

Case (vi): The estimates in this case follow from Lemma~\ref{lem:-solsm} and Lemma~\ref{lem:Nbds}.
\end{Proof}

\subsection{ The components $\bG_\lambda^{12}(x,y)$ and $\bG_\lambda^{22}(x,y)$ }
The analysis for the components $\bG_\lambda^{12}(x,y)$ and $\bG_\lambda^{22}(x,y)$  of $\bG_\lambda(x,y)$ proceeds along similar lines as in Section~\ref{sec:G11}.  The formulas for the coefficients in this case are,
\[  \left(\begin{array}{c} c_1 \\ c_2 \\ c_3 \\ c_4 \end{array}\right)= \left(\begin{matrix} \phi_1^+(y) &  \phi_2^+(y) &\phi_1^-(y) &\phi_2^-(y)\end{matrix} \right)^{-1} \left(\begin{array}{c} 0 \\ 0 \\ 0 \\ -\frac{1}{\sigma} \end{array}\right),\]
from which we obtain
\begin{align*} 
c_1(y)&= \frac{1}{\sigma\W_\lambda(y)}\mathrm{det}\left(\begin{matrix} \pi_4 \phi_2^+(y) & \pi_4 \phi_1^-(y) & \pi_4\phi_2^-(y)\end{matrix}\right), \\
c_2(y)&= -\frac{1}{\sigma\W_\lambda(y)}\mathrm{det}\left(\begin{matrix} \pi_4 \phi_1^+(y) & \pi_4 \phi_1^-(y) & \pi_4\phi_2^-(y)\end{matrix}\right), \\
c_3(y)&= \frac{1}{\sigma\W_\lambda(y)}\mathrm{det}\left(\begin{matrix} \pi_4 \phi_1^+(y) & \pi_4 \phi_2^+(y) & \pi_4\phi_2^-(y)\end{matrix}\right) ,\\
c_4(y)&= -\frac{1}{\sigma\W_\lambda(y)}\mathrm{det}\left(\begin{matrix} \pi_4 \phi_1^+(y) & \pi_4 \phi_2^+(y) & \pi_4\phi_1^-(y)\end{matrix}\right).
\end{align*}
From this we determine 
\begin{eqnarray} \bG_\lambda^{12}(x,y)&=&\left\{ \begin{array}{cc} \frac{M_1^+(y)\phi_{1,u}^+(x)}{\sigma\W_\lambda(y) }+ \frac{M_2^+(y)\phi_{2,u}^+(x)}{\sigma\W_\lambda(y) },   & x>y, \\
\frac{M_1^-(y)\phi_{1,u}^-(x)}{\sigma\W_\lambda(y) }+ \frac{M_2^-(y)\phi_{2,u}^-(x)}{\sigma\W_\lambda(y) },   & x<y, \end{array}\right. \label{eq:G12} \\
 \bG_\lambda^{22}(x,y)&=& \left\{ \begin{array}{cc} \frac{M_1^+(y)\phi_{1,v}^+(x)}{\sigma\W_\lambda(y) }+ \frac{M_2^+(y)\phi_{2,v}^+(x)}{\sigma\W_\lambda(y) },   & x>y, \\
\frac{M_1^-(y)\phi_{1,v}^-(x)}{\sigma\W_\lambda(y) }+ \frac{M_2^-(y)\phi_{2,v}^-(x)}{\sigma\W_\lambda(y) } ,  & x<y .\end{array}\right. \label{eq:G22} \end{eqnarray}
In both of these expressions,
\bqq
\begin{split} M_1^+(y)&=\mathrm{det}\left(\begin{matrix} \pi_4 \phi_2^+(y) & \pi_4 \phi_1^-(y) & \pi_4\phi_2^-(y)\end{matrix}\right), \\
M_2^+(y)&=- \mathrm{det}\left(\begin{matrix} \pi_4 \phi_1^+(y) & \pi_4 \phi_1^-(y)  \pi_4\phi_2^-(y)\end{matrix}\right), \\
M_1^-(y)&= -\mathrm{det}\left(\begin{matrix} \pi_4 \phi_1^+(y) & \pi_4 \phi_2^+(y) & \pi_4\phi_2^-(y)\end{matrix}\right), \\
M_2^-(y)&=  \mathrm{det}\left(\begin{matrix} \pi_4 \phi_1^+(y) & \pi_4 \phi_2^+(y) & \pi_4\phi_1^-(y)\end{matrix}\right).
\end{split}
\label{eqMij}
\eqq
\begin{lem}\label{lem:Mbds}Assume that $\lambda$ is to the right of $\Gamma$ off the imaginary axis and with $|\lambda|<M_s$. Then, the terms $M_{j}^\pm(y)$ from \eqref{eqMij} obey the same bounds as the $N_j^\pm(y)$ from Lemma~\ref{lem:Nbds} for $j=1,2$.  In addition, $\bG_\lambda^{12}(x,y)$ and $\bG_\lambda^{22}(x,y)$ also satisfy the same bounds as $\bG_\lambda^{11}(x,y)$ and $\bG_\lambda^{12}(x,y)$ in Lemma~\ref{lem:G11bds}.
\end{lem}
\begin{Proof} Note that  $M_j^\pm(y)$ are simply $N_j^\pm(y)$ with $\pi_2$ replaced with $\pi_4$.  Since the proof of Lemma~\ref{lem:Nbds} and Lemma~\ref{lem:G11bds} did not rely on any properties of $\pi_2$, the same bounds hold here.
\end{Proof}

\subsection{Proof of Propostion~\ref{thmGlambda}}

To conclude the proof of Proposition~\ref{thmGlambda}, we need only to establish estimate (\ref{eqGlambdanearzero}) for those values of $\lambda$ to the right of $\Gamma$, off the imaginary axis and within the ball $|\lambda|<M_s$. We observe that the result follows directly from Lemma~\ref{lem:G11bds} and Lemma~\ref{lem:Mbds}.  Note that the estimate provided in (\ref{eqGlambdanearzero}) is not the sharpest possible, but will be sufficient to establish our main result.

\section{ The temporal Green's function $\cG(t,x,y)$. }
In this section, we derive bounds for the the temporal Green's function 
\[ \cG(t,x,y)=\left(\begin{array}{cc} \cG^{11}(t,x,y) & \cG^{12} (t,x,y) \\
\cG^{21}(t,x,y) & \cG^{22}(t,x,y) \end{array}\right), \quad t>0, \quad x,y\in\R. \]
We recall that $\cG(t,x,y)$ can be recovered from the inverse Laplace transform of $\bG_\lambda(x,y)$,
\[ \cG(t,x,y)=\frac{1}{2\pi \mathbf{i}} \int_\mathscr{C} e^{\lambda t}\bG_\lambda(x,y)
\mathrm{d}\lambda,\]
for some well-chosen contour $\mathscr{C}$ which does not intersect with the spectrum of $\cL$.

The goal is to use the bounds in Proposition~\ref{thmGlambda} in combination with advantageous choice for the integration contour $\mathscr{C}$ to obtain sufficient bounds on the temporal Green's function so as to perform a nonlinear stability argument in the following section.  The main result of this section is the following Proposition.  

\begin{prop}\label{timeGreen}
Under the assumptions of our main theorem, and for some constants $\kappa>0$, $r>0$ and $C>0$, the Green's function $\cG(t,x,y)$ for $\partial_t p = \cL p$ satisfies the following estimates.
\begin{itemize}
\item[(i)] For $|x-y|\geq Kt$ or $t<1$,  with $K$ sufficiently large, 
\bqs
\left|\cG^{ij}(t,x,y)\right|\leq C\frac{1}{t^{1/2}}e^{-\frac{|x-y|^2}{\kappa t}},
\eqs
for all $i,j\in\left\{1,2\right\}$.
\item[(ii)] For $|x-y|\leq Kt$ and $t\geq 1$, with $K$ as above, 
\bqs
\left|\cG^{ij}(t,x,y)\right|\leq C\left(\frac{1+|x-y|}{t^{3/2}} \right)e^{-\frac{|x-y|^2}{\kappa t}}+Ce^{-rt},
\eqs
for all $i,j\in\left\{1,2\right\}$.
\end{itemize}
\end{prop}

\begin{Proof} 
Case (i): Consider first the case when $|x-y|\geq Kt$ or $t<1$.  The proof of this case is standard and we refer the reader to the Proof of Theorem 8.3 in  \cite{zumbrun98} or Proposition 4.1 of \cite{faye19} for a proof that can be modified to the current context in a straightforward fashion.

Case (ii): Now consider the scenario where $|x-y|\leq Kt$ and $t\geq 1$.   The proof in this case is similar to that of Proposition 4.1 in  \cite{faye19}, which in turn mimics the approach taken in a number of previous works; see once again \cite{zumbrun98}.  We include some details here for completeness and for the convenience of the reader. 

Recall the definition of the region
\[ \Omega_\delta=\left\{ \lambda\in\mathbb{C}\ \  | \  \mathrm{Re}(\lambda)\geq -\delta_0-\delta_1 |\mathrm{Im}(\lambda)| \right\} , \]
with $\delta_0$ and $\delta_1$ chosen as in Lemma~\ref{lemnozeros} together with its boundary $\Gamma_\delta$.  We will perform the Laplace inversion on a contour $\mathscr{C}$ which is given by a parabolic segment near the origin followed by the linear contour $\Gamma_\delta$.  This contour will be divided into four segments $\mathscr{C}=\Gamma_1\cup\Gamma_2\cup\Gamma_3\cup\Gamma_4$.
\begin{itemize}
\item $\Gamma_1$ consists of a parabolic contour defined as follows.  Let $\rho>0$ and consider 
\[ \sqrt{\lambda}=\rho+\rmi k, \]
so that 
\[ \lambda=\rho^2-k^2+2\rmi \rho k. \]
Since we are interested in bounds for which  $|x-y|<Kt$, we will take 
\[ \rho=\frac{|x-y|}{Lt}\]
with $L$ chosen sufficiently large so that $\Gamma_1$ is contained in the region where the "small  $\lambda$" estimates of Proposition~\ref{thmGlambda} hold and to the right of $\Gamma_\delta$. 
\item $\Gamma_2$ consists of a linear contour following $\Gamma_\delta$.  This contour is defined for those values of $\lambda$ where the "small  $\lambda$" estimates in part (i) of  Proposition~\ref{thmGlambda} hold. 
\item $\Gamma_3$ is the continuation along the contour $\Gamma_\delta$ where the "medium $\lambda$" estimates detailed in part (iii) of Proposition~\ref{thmGlambda} are valid.
\item $\Gamma_4$ is the continuation along the contour $\Gamma_\delta$ where the "large $\lambda$" estimates detailed in part (ii) of Proposition~\ref{thmGlambda} hold.  
\end{itemize}
We begin with the estimates along $\Gamma_1$.  We first note that
\bqs
\md \lambda = 2\mbi \left(\frac{(x-y)}{Lt}+\mbi k\right)\md k \quad  \text{ and } \quad \lambda = \frac{(x-y)^2}{L^2t^2}-k^2+2\frac{(x-y)}{Lt}\mbi k.
\eqs
Then
\[ \frac{1}{2\pi \mbi }\int_{\Gamma_1} e^{\lambda t} \bG_\lambda(x,y)\mathrm{d}\lambda  =  
 \frac{1}{\pi  }e^{\rho^2 t-\rho(x-y)} \int_{-k^*}^{k^*} e^{-k^2 t}e^{2\mbi \rho k -\mbi k(x-y)} \boldH_{\lambda(k)}(x,y)(\rho+\mbi k) \mathrm{d}k, \]
from which we expand $\boldH_{\lambda(k)}$ into its real and imaginary parts, $\boldH_{\lambda(k)}:=H_R(x,y,k)+\rmi H_I(x,y,k)$ so that the integral reduces to 
\[ \frac{1}{\pi  }e^{\rho^2 t-\rho(x-y)} \int_{-k^*}^{k^*} e^{-k^2 t} \left(H_R(x,y, k)+\mbi H_I(x,y,k)\right)(\rho+\mbi k) \mathrm{d}k. \]
Here $k^*$ prescribes the end of the contour $\Gamma_1$ and the start of contour $\Gamma_2$.  We compute its value below, but taking it to be arbitrary here does not change the analysis.  

Since  $H_R$ is bounded we have that 
\[ \left|\frac{1}{\pi  }e^{\rho^2 t-\rho(x-y)} \int_{-k^*}^{k^*} e^{-k^2 t} \rho H_R^{ij}(x,y, k) \mathrm{d}k\right|\leq C\frac{\rho}{\sqrt{t}} e^{\rho^2 t-\rho(x-y)},  \]
for all $i,j\in\left\{1,2\right\}$. On the other hand, since $H_I(x,y,k)$ is odd in $k$, it can be expressed as $H_I(x,y,k)=k\widetilde{H}_I(x,y,k)$, where $\widetilde{H}_I(x,y,k)$ is again bounded.  Therefore,  
\[ \left|\frac{1}{\pi  }e^{\rho^2 t-\rho(x-y)} \int_{-k^*}^{k^*} e^{-k^2 t} k^2 \widetilde{H}_I^{ij}(x,y, k) \mathrm{d}k\right| \leq \frac{C}{t^{3/2}} e^{\rho^2 t-\rho(x-y)},\]
for all $i,j\in\left\{1,2\right\}$.
Recall that  $\rho=\frac{|x-y|}{Lt}$, so that we have obtained the estimate for all $i,j\in\left\{1,2\right\}$
\[  \frac{1}{2\pi }\left|\int_{\Gamma_1} e^{\lambda t} \bG_\lambda^{ij}(x,y)\mathrm{d}\lambda\right| \leq C\frac{1+|x-y|}{t^{3/2}} e^{-\frac{|x-y|^2}{\kappa t}},\]
for $\kappa=L^2/(L-1)$.  

We now turn to the analysis of the contour integral along $\Gamma_2$. Since the analysis is equivalent, we focus only on the segment in the positive half plane.  There, the contour $\Gamma_2$ can be parameterized by 
\[ \lambda=-\delta_0 +\cos(\theta) \ell+\mbi \sin(\theta)\ell, \]
for some fixed $\pi/2<\theta<\pi$.  The contours $\Gamma_1$ and $\Gamma_2$ intersect at 
\bqq k^*(\rho)=-\rho \cot(\theta)+\sqrt{\rho^2 \csc^2(\theta)+\delta_0},\ \quad \ell_1(\rho)=\frac{2\rho k^*(\rho)}{\sin(\theta)}.\label{eql1} \eqq
We now estimate for each $i,j\in\left\{1,2\right\}$
\begin{align} \frac{1}{2\pi}\left|\int_{\Gamma_2} e^{\lambda t} \bG_\lambda^{ij}(x,y)\mathrm{d}\lambda\right|  &\leq   
 Ce^{-\delta_0 t} \int_{\ell_1}^{\ell_2} e^{t\cos(\theta)\ell -(x-y)\mathrm{Re}(\sqrt{-\delta_0+\cos(\theta)\ell+\mbi \sin(\theta)\ell}) }\mathrm{d}\ell \nonumber \\
&\leq Ce^{-\delta_0 t} \int_{\ell_1}^{\ell_2} e^{t\cos(\theta)\ell } \mathrm{d}\ell \leq C \frac{e^{-\delta_0 t}}{t}e^{t\cos(\theta)\ell_1}  ,
\end{align}
we were have used that $(x-y)\mathrm{Re}(\sqrt{-\delta_0+\cos(\theta)\ell+\mbi \sin(\theta)\ell}) >0$ and integrated noting that $\cos(\theta)<0$.    Now, by virtue of (\ref{eql1}) we observe that 
\[ t\cos(\theta) \ell_1(\rho) < 2t\cot(\theta)\rho^2\left( -\cot(\theta)+\csc(\theta)\right)<0.\]
We therefore obtain the bound
\[ \frac{1}{2\pi}\left|\int_{\Gamma_2} e^{\lambda t} \bG_\lambda^{ij}(x,y)\mathrm{d}\lambda\right|  \leq C\frac{e^{-\delta_0 t}}{t} e^{-\frac{|x-y|^2}{\kappa t}} \leq C\frac{1}{t^{3/2}} e^{-\frac{|x-y|^2}{\kappa t}},\]
for some $\kappa>0$ and all $i,j\in\left\{1,2\right\}$.

Next consider the integral along $\Gamma_3$.  Here the "medium" $\lambda$ estimate holds and we simply have that $\bG_\lambda(x,y)$ is uniformly bounded in this region.  We then calculate, 
\[ \frac{1}{2\pi }\left|\int_{\Gamma_3} e^{\lambda t} \bG_\lambda^{ij}(x,y)\mathrm{d}\lambda\right| \leq Ce^{-rt}, \]  
for each $i,j\in\left\{1,2\right\}$ and for some $r>0$.  

Finally, we consider the integral along $\Gamma_4$.   The large $\lambda$ bounds apply here and the analysis follows as in the case of (i) above.  We find that each component satisfies
\[ \frac{1}{2\pi }\left|\int_{\Gamma_4} e^{\lambda t} \bG_\lambda^{ij}(x,y)\mathrm{d}\lambda\right| \leq C e^{-\delta_0 t} \int_{\ell_3}^\infty \frac{1}{\sqrt{\ell}} e^{\cos(\theta)\ell t}\mathrm{d}\ell\leq C e^{-rt} , \]  
for some $r>0$.  
\end{Proof}

\section{Nonlinear Stability}
We now turn to the question of nonlinear stability and establish Theorem~\ref{thmmain}.  We consider solutions of the nonlinear system of equations (\ref{eq:pq}). The fact that $\cL$ generates an analytic semigroup implies that the Cauchy problem associated to (\ref{eq:pq}) with initial condition $(p_0,q_0)\in L^1(\R)\cap L^\infty(\R)$ with $\int_\R |y|\left(|p_0(y)| + |q_0(y)|\right)\mathrm{d}y<+\infty$ is locally well-posed in $L^\infty(\R)$. We let $T^*>0$ be the associated maximal time of existence of a such a solution $(p,q)$.

We begin by formulating (\ref{eq:pq}) in its integral form for all $t\in(0,T^*)$ and $x\in\R$
\bqq
\left(\begin{matrix}p(t,x) \\ q(t,x) \end{matrix}\right)  = \int_\R \cG(t,x,y)\left(\begin{matrix}p_0(y) \\ q_0(y) \end{matrix}\right)\mathrm{d}y+\int_0^t \int_\R \cG(t-\tau,x,y)\left( \begin{matrix}  \mathcal{N}_u( p(\tau,y),q(\tau,y)) \\  \mathcal{N}_v(p(\tau,y),q(\tau,y))\end{matrix}\right) \mathrm{d}y\mathrm{d}\tau.
\label{inteqp}
\eqq
We will require the following estimates regarding the Green's function $\cG(t,x,y)$.
\begin{lem}
For  $t<1$ and all $x\in\R$, we have that
\bqq 
\left| \int_\R \cG^{ij}(t,x,y)h(y)\mathrm{d}y \right| \leq C \| h\|_{L^\infty(\R)}, \quad i,j\in\left\{1,2\right\}.
\label{eqGestsmall}
\eqq
Conversely, for $t\geq 1$ and all $x\in\R$, we have that
\bqq \left| \int_\R \cG^{ij}(t,x,y)h(y)\mathrm{d}y \right| \leq C \frac{1+|x|}{(1+t)^{3/2}}\int_\R(1+|y|)|h(y)|\md y, \quad i,j\in\left\{1,2\right\}.\label{eqGest} \eqq
\end{lem}

\begin{Proof}
Both estimates are a direct consequence of Proposition~\ref{timeGreen}.  For the second estimate, we have
\begin{align}
\left| \int_\R \cG^{ij}(t,x,y)h(y)\mathrm{d}y \right| &\leq \int_{-\infty}^{x-Kt} C\frac{1}{t^{1/2}}e^{-\frac{|x-y|^2}{\kappa t}}|h(y)|\mathrm{d}y \nonumber \\
&+\int_{x-Kt}^{x+Kt} C\left( \left(\frac{1+|x-y|}{t^{3/2}} \right)e^{-\frac{|x-y|^2}{\kappa t}}+e^{-rt}\right)|h(y)|\mathrm{d}y  \nonumber \\
&+\int_{x+Kt}^{\infty }C\frac{1}{t^{1/2}}e^{-\frac{|x-y|^2}{\kappa t}}|h(y)|\mathrm{d}y.
\end{align}
The first and last integral decay exponentially in time, while for the middle integral 
we take the limits of integration to infinity and note that $1+|x-y|\leq 1+|x|+|y|+|x||y|$ from which we obtain
\[ \left| \int_\R \cG^{ij}(t,x,y)h(y)\mathrm{d}y \right| \leq Ce^{-rt}\int_{\mathbb{R}}(1+|y|)|h(y)|\mathrm{d}y +C \frac{1+|x|}{(1+t)^{3/2}}\int_\R(1+|y|)|h(y)|\md y.\]
The exponential decay of the first integral allows it to be incorporated into the second and we obtain our desired estimate.  
\end{Proof}
We now return to (\ref{inteqp}) and for $t\in [0,T^*)$  we define
\bqs
\Theta(t) := \underset{0\leq \tau \leq t}{\sup}~\underset{x\in\R}{\sup} ~ \frac{ (1+\tau )^{3/2}}{1+|x|}\left( |p(\tau,x)| + |q(\tau,x)|\right).
\eqs
We consider first the case of $0<t<1$.  Here we apply the estimate (\ref{eqGestsmall}) to (\ref{inteqp}).  For the first integral it holds that 
\[ \left\| \int_\R \cG(t,x,y)\left(\begin{array}{c}p_0(y) \\ q_0(y) \end{array}\right)\mathrm{d}y \right\| \leq C \left(\|p_0\|_{L^\infty(\R)} + \|q_0\|_{L^\infty(\R)} \right). \]
For the second, we have 
\begin{eqnarray*}  \left\| \int_\R \cG(t-\tau,x,y)\left( \begin{array}{c}  - p(\tau,y)( \omega(y) p(\tau,y)+a\omega(y) q(\tau,y)) \\  -rq(\tau,y)(b\omega(y) p(\tau,y)+\omega(y) q(\tau,y))\end{array}\right) \mathrm{d}y \right\|\\
\leq C \frac{\Theta(t)^2}{(t-\tau)^{1/2}} \int_\R e^{-\frac{|x-y|^2}{\kappa(t-\tau)}} (1+|y|)^2 \omega(y) \mathrm{d}y .
\end{eqnarray*}
Now integrating we obtain, again for $0<t<1$, 
\begin{eqnarray*}  \left\|\int_0^t \int_\R \cG(t-\tau,x,y)\left( \begin{array}{c}  - p(\tau,y)( \omega(y) p(\tau,y)+a\omega(y) q(\tau,y)) \\  -rq(\tau,y)(b\omega(y) p(\tau,y)+\omega(y) q(\tau,y))\end{array}\right) \mathrm{d}y \mathrm{d}\tau\right\| \\
\leq C \Theta(t)^2 \int_\R  (1+|y|)^2 \omega(y) \mathrm{d}y .
\end{eqnarray*}
Since the weight $\omega(y)$ is exponentially localized, we obtain the following inequality, valid for $t<1$
\bqq
 \Theta(t)\leq C \left(\|p_0\|_{L^\infty(\R)} + \|q_0\|_{L^\infty(\R)} \right) +C  \Theta(t)^2. 
\label{eqt<1}
 \eqq
We now consider the case of $t\geq 1$.  Focusing on the first integral in (\ref{inteqp}), we apply  (\ref{eqGest}) and obtain 
\[ \left\|\int_\R \cG(t,x,y)\left(\begin{array}{c}p_0(y) \\ q_0(y) \end{array}\right)\mathrm{d}y \right\| \leq C\frac{1+|x|}{(1+t)^{3/2}}\int_\R (1+|y|)\left(|p_0(y)| + |q_0(y)|\right)\mathrm{d}y.  \]
For the second integral, we have 
\begin{eqnarray*}  \left\| \int_\R \cG(t-\tau,x,y)\left( \begin{array}{c}  - p(\tau,y)( \omega(y) p(\tau,y)+a\omega(y) q(\tau,y)) \\  -rq(\tau,y)(b\omega(y) p(\tau,y)+\omega(y) q(\tau,y))\end{array}\right) \mathrm{d}y \right\|\\
\leq C \frac{\Theta(t)^2}{(1+t-\tau)^{3/2}}\frac{1+|x|}{(1+\tau)^3} \int_\R (1+|y|)^3 \omega(y) \mathrm{d}y .
\end{eqnarray*}
Integrating with respect to $\tau$ we note that (see for instance \cite{xin92}) 
\bqs
\int_0^t \frac{1}{(1+t-\tau)^{3/2}(1+\tau)^3}\md \tau \leq \frac{\tilde{C}}{(1+t)^{3/2}}.
\eqs
This gives the inequality, valid for $t\geq 1$,
\bqq 
\Theta(t)\leq C\int_\R (1+|y|)\left(|p_0(y)| + |q_0(y)|\right)\mathrm{d}y +C\Theta(t)^2\int_\R (1+|y|)^3 \omega(y) \mathrm{d}y. 
\label{eqt>1}
\eqq
Combining \eqref{eqt<1} and \eqref{eqt>1}, we get that there exist $C_1>0$ and $C_2$ such that for all $t\in [0,T^*)$
\bqq
\Theta(t)\leq C_1 \Omega +C_2 \Theta(t)^2,
\label{eqallt}
\eqq
with
\bqs
\Omega:=\left\| (p_0,q_0)  \right\|_{L^\infty(\R)}+\left\| (1+|\cdot|) (p_0,q_0)  \right\|_{L^1(\R)}.
\eqs
As a consequence, if we assume that the initial perturbation $(p_0,q_0)$ is small enough so that
\bqs
2C_1\Omega < 1, \quad \text{ and } \quad 4C_1C_2 \Omega <1,
\eqs
then we claim that $\Theta(t)\leq 2C_1 \Omega<1$ for all $t\in [0,T^*)$. This implies that the maximal time of existence is $T^*=\infty$ and the solution $(p,q)$ of (\ref{eq:pq}) satisfies
\bqs
\underset{t\geq 0}{\sup}~\underset{x\in\R}{\sup} ~ \frac{ (1+\tau )^{3/2}}{1+|x|}\left( |p(\tau,x)| + |q(\tau,x)|\right)<2C_1\Omega,
\eqs
which concludes the proof of Theorem~\ref{thmmain}. Coming back to the claim, we see that by eventually taking $C_1$ even larger, we can always assume that at $t=0$
\bqs
\Theta(0)=\underset{x\in\R}{\sup} ~ \frac{ |p_0(x)|+|q_0(x)|}{1+|x|} \leq \left\| (p_0,q_0)  \right\|_{L^\infty(\R)} <\Omega < 2C_1\Omega,
\eqs
such that by continuity of $\Theta(t)$ we have that for small time $\Theta(t)<2C_1\Omega$. Suppose there exists $T>0$ where $\Theta(T)=2C_1\Omega$ for the first time, then from \eqref{eqallt} we have
\bqs
\Theta(T)\leq C_1\Omega(1+4C_1C_2\Omega)<2C_1\Omega,
\eqs
which gives a contradiction and proves the claim.

\section*{Acknowledgments} 
GF received support from the ANR project NONLOCAL ANR-14-CE25-0013. MH received partial support from the National Science Foundation through grant NSF-DMS-1516155.

\appendix

\section{Asymptotic behavior of the critical traveling front $(U_*,V_*)$}

In this section, we collect some known results regarding the precise asymptotic behavior of the critical traveling front $(U_*,V_*)$ solution of \eqref{eq:wave}-\eqref{eq_limits} that can be found in \cite{kanon97,morita09,girardin18}. To this aim, we introduce four different dispersion relations
\begin{align*}
d_u^{-\infty}(\lambda)&=\lambda^2+c_*\lambda-1,\\
d_v^{-\infty}(\lambda)&=\sigma\lambda^2+c_*\lambda+r(1-b),\\
d_u^{+\infty}(\lambda)&=\lambda^2+c_*\lambda+(1-a),\\
d_v^{+\infty}(\lambda)&=\sigma\lambda^2+c_*\lambda-r.
\end{align*}
We define 
\bqs
\mu_u^{-\infty}=-\gamma_*+\sqrt{\gamma_*^2+1}>0 \quad \text{ and } \quad \mu_v^{-\infty}=-\frac{\gamma_*}{\sigma}+\frac{1}{\sigma}\sqrt{\gamma_*^2+\sigma r(b-1)}>0,
\eqs
together with
\bqs
\nu_u^{+\infty}=-\gamma_*<0 \quad \text{ and } \quad \nu_v^{+\infty}=-\frac{\gamma_*}{\sigma}-\frac{1}{\sigma}\sqrt{\gamma_*^2+\sigma r}<0.
\eqs

\paragraph{Asymptotic behavior at $+\infty$.} Regarding the asymptotic decay at $+\infty$, one needs to compare $\nu_u^{+\infty}$ and $\nu_v^{+\infty}$. From assumption \eqref{par_lin} and the fact that $0<\sigma<2$, we always have
\bqs
(\sigma-1)\gamma_*<\sqrt{\gamma_*^2+\sigma r},
\eqs
which in turn implies that $\nu_v^{+\infty}<\nu_u^{+\infty}$ which falls into case (2)-(b)-(iii) of \cite[Lemma A.1]{girardin18}. As a consequence, we have
\bqq \left( \begin{array}{c} U_*(\xi) \\ V_*(\xi)-1 \end{array}\right) = \left( \begin{array}{c} \beta  \\ \dfrac{rb \beta}{d_v^{+\infty}(\nu_u^{+\infty})} \end{array}\right)\xi e^{-\gamma_* \xi} +\mathrm{o}\left(\xi e^{-\gamma_* \xi}\right) \quad \text{as $\xi\to+\infty$}, 
\label{eq:asympt_front+}
\eqq
for some $\beta>0$. Let us remark that $d_v^{+\infty}(\nu_u^{+\infty})=(\sigma-2)\gamma_*^2-r<0$.

\paragraph{Asymptotic behavior at $-\infty$.} Regarding the asymptotic behavior at $-\infty$ we have the following classification. There exist $\beta_1>0$ and $\beta_2>0$ such that, as $\xi\rightarrow-\infty$
\begin{itemize}
\item if $\mu_u^{-\infty}<\mu_v^{-\infty}$, then
\bqq 
\left( \begin{array}{c} 1-U_*(\xi) \\ V_*(\xi) \end{array}\right)= \left( \begin{array}{c} \beta_1e^{\mu_u^{-\infty}\xi}  \\ \beta_2 e^{\mu_v^{-\infty}\xi} \end{array}\right)+\mathrm{h.o.t.}~;
\label{eq:asympt_front-1}
\eqq
\item if $\mu_u^{-\infty}>\mu_v^{-\infty}$, then $d_u^{-\infty}(\mu_v^{-\infty})<0$ and 
\bqq 
\left( \begin{array}{c} 1-U_*(\xi) \\ V_*(\xi) \end{array}\right)= \left( \begin{array}{c} -\dfrac{a\beta_2}{d_u^{-\infty}(\mu_v^{-\infty})}  \\ \beta_2 \end{array}\right)e^{\mu_v^{-\infty}\xi}+\mathrm{h.o.t.}~;
\label{eq:asympt_front-2}
\eqq
\item if $\mu_u^{-\infty}=\mu_v^{-\infty}$, then 
\bqq 
\left( \begin{array}{c} 1-U_*(\xi) \\ V_*(\xi) \end{array}\right)= \left( \begin{array}{c} -\beta_2 \xi   \\ \dfrac{2\sqrt{\gamma_*^2+1}\beta_2}{a} \end{array}\right)e^{\mu_v^{-\infty}\xi}+\mathrm{h.o.t.}~.
\label{eq:asympt_front-3}
\eqq
\end{itemize}

\section{Spectral properties of $\tcL$ and $\cL$}

In this section, we investigate the spectral properties of the operator $\tcL$ defined as the linearization of the Lotka-Volterra system around the critical traveling front $(U_*,V_*)$ 
\bqs
\widetilde{\mathcal{L}}=\left( \begin{array}{cc} \partial_{xx}+c_*\partial_x+(1-2U_*(x)-aV_*(x)) & -aU_*(x) \\
-rbV_*(x) & \sigma\partial_{xx}+c_*\partial_x+r(1-bU_*(x)-2V_*(x)) \end{array}\right),
\eqs
both set on $L^2(\R)\times L^2(\R)$ with dense domain $H^2(\R)\times H^2(\R)$ or on the weighted Sobolev space $L^2_{\omega^{-1}} (\R)\times L^2_{\omega^{-1}}(\R)$ with dense domain $H^2_{\omega^{-1}}(\R)\times H^2_{\omega^{-1}}(\R)$  where
\bqs
L^2_{\omega^{-1}} (\R)=\left\{ f \in L^2(\R) ~|~ (\omega^{-1}f) \in L^2(\R) \right\},
\eqs
and
\bqs
H^2_{\omega^{-1}} (\R)=\left\{ f \in L^2(\R) ~|~ f,f',f'' \in L^2_{\omega^{-1}} (\R) \right\},
\eqs
with the weight $\omega>0$ from \eqref{eq:weights}. The first result asserts that under our assumptions on the coefficients $\lambda=0$ is not an eigenvalue of $\tcL$ in the exponentially weighted space $L^2_{\omega^{-1}} (\R)\times L^2_{\omega^{-1}}(\R)$.

\begin{lem}\label{lem:0eig}
Assume that assumptions \eqref{par_mono}-\eqref{par_lin} hold. Then $\lambda=0$ is not an eigenvalue of $\tcL$ set on $L^2_{\omega^{-1}} (\R)\times L^2_{\omega^{-1}}(\R)$.
\end{lem}

The second result excludes the possibility of having eigenvalues with positive real part for $\tcL$ set on $L^2_{\omega^{-1}} (\R)\times L^2_{\omega^{-1}}(\R)$.  In fact, any eigenvalues of $\tcL$ either on $L^2(\R)\times L^2(\R)$ or on the weighted Sobolev space $L^2_{\omega^{-1}} (\R)\times L^2_{\omega^{-1}}(\R)$ has to be real. This is a consequence of the monotonicity property of $\tcL$, namely that the off-diagonal terms of $\tcL$ are strictly negative (whether the operator is seen on $L^2(\R)\times L^2(\R)$ or $L^2_{\omega^{-1}} (\R)\times L^2_{\omega^{-1}}(\R)$). We also introduce the conjugate operator 
\bqs
\cL = \omega^{-1} \tcL \omega, \quad \mathcal{D}(\cL)=H^2(\R)\times H^2(\R).
\eqs
We have the following result. 

\begin{lem}\label{lem:noposeig}
Let $\lambda\in\R$ be a real eigenvalue of $\cL$ defined on $L^2(\R)\times L^2(\R)$, then $\lambda<0$.
\end{lem}

Note that the proofs of Lemma~\ref{lem:0eig}-\ref{lem:noposeig} are very close to the ones that can be found in \cite{bates2006spectral,leung11}. For the sake of clarity and completeness, we have reproduced them here.

Before proceeding to the proof of the two lemmas, we first change coordinates such that system \eqref{eq:main} is of cooperative type. Thus, we let $w=1-v$ and we find
\bqq
\begin{split}
u_t&= u_{xx}+u(1-a-u+aw),\\
w_t&=  \sigma w_{xx}+r(1-w)(bu-w).
\end{split}
\label{eq:mainco} 
\eqq
Next, we note that the critical front $(U_*,V_*)$ is transformed into $(U_*,W_*)$ and connects $(1,1)$ at $-\infty$ to $(0,0)$ at $+\infty$ with asymptotics at $+\infty$ given by \eqref{eq:asympt_front+}
\bqs \left( \begin{array}{c} U_*(\xi) \\ W_*(\xi) \end{array}\right) = \left( \begin{array}{c} \beta  \\ -\dfrac{rb \beta}{d_v^{+\infty}(\nu_u^{+\infty})} \end{array}\right)\xi e^{-\gamma_* \xi} +\mathrm{o}\left(\xi e^{-\gamma_* \xi}\right) \quad \text{as $\xi\to+\infty$}, 
\eqs
for some $\beta>0$ and $d_v^{+\infty}(\nu_u^{+\infty})<0$. Let also note that both $U_*'<0$ and $W_*'<0$, while $0<U_*,W_*<1$. We denote by $\tscL$ the linearized operator around the critical front $(U_*,W_*)$
\bqs
\tscL:=\left( \begin{array}{cc} \partial_{xx}+c_*\partial_x+(1-a-2U_*(x)+aW_*(x)) & aU_*(x) \\
rb(1-W_*(x)) & \sigma\partial_{xx}+c_*\partial_x+r(-1-bU_*(x)+2W_*(x)) \end{array}\right),
\eqs
and we note that $\tscL(U_*',W_*')=0$. We will consider $\tscL$ either on $L^2(\R)\times L^2(\R)$ or on the weighted Sobolev space $L^2_{\omega^{-1}} (\R)\times L^2_{\omega^{-1}}(\R)$ with dense domains $H^2(\R)\times H^2(\R)$ or $H^2_{\omega^{-1}} (\R)\times H^2_{\omega^{-1}}(\R)$ respectively. 

Throughout the end of this section, we will use the following ordering
\begin{align*}
(p_1,q_1) &\succ (p_2,q_2) \Longleftrightarrow p_1(x) > p_2(x) \quad \text{ and } \quad q_1(x) > q_2(x) \quad \text{for all} \quad x\in\R,\\
(p_1,q_1) &\succeq (p_2,q_2) \Longleftrightarrow p_1(x) \geq p_2(x) \quad \text{ and } \quad q_1(x) \geq q_2(x) \quad \text{for all} \quad x\in\R.
\end{align*}

\subsection{Proof of Lemma~\ref{lem:0eig}}

Suppose there is a nonzero function $\cV\in H^2_{\omega^{-1}} (\R)\times H^2_{\omega^{-1}}(\R)$ which is a solution of $\tscL \cV = 0$. The strategy is to compare $\cV$ to $\U_*=(-U_*',-W_*')  \succ (0,0)$ which also satisfies $\tscL \U_*=0$. More precisely, we will show that if such a $\cV=(v_1,v_2)\in H^2_{\omega^{-1}} (\R)\times H^2_{\omega^{-1}}(\R)$ exists then necessarily we have 
\bqs
|\tau\cV(x)|:=(|\tau v_1(x)|,|\tau v_2(x)|) \preceq \U_*(x)
\eqs
 for all $\tau \in\R$ and all $x\in\R$ which is obviously a contradiction to the fact that $\U_*$ is bounded. We introduce the following set
\bqs
\mathcal{S}=\left\{ \tau\in\R ~|~ |\tau\cV(x)| \preceq \U_*(x), \quad x\in\R \right\}.
\eqs
As $0\prec \U_*$ we have that $0\in \mathcal{S}$ and so $\mathcal{S}$ is non empty. It is also straightforward to check that it is closed. Thus, we are going to show that is open. Let $\tau \in \mathcal{S}$ be fixed.

\textbf{Step 1.} We have $|\tau\cV(x)| \preceq \U_*(x)$ for all $x\in\R$. We claim that it implies that $|\tau\cV(x)| \prec \U_*(x)$ for all $x\in\R$. We consider $\cW(x)=\U_*(x)-\tau\cV(x)\succeq (0,0)$ which satisfies $\tscL \cW =0$ with $\cW(-\infty)=\cW(+\infty)=(0,0)$. Noticing that $aU_*>0$ and $rb(1-W_*)>0$, the maximum principle implies that $\cW\succ (0,0)$ unless it is identically zero which is impossible because in that case one would have $|\tau\cV(x)| = \U_*(x)$ but $\U_*\notin H^2_{\omega^{-1}} (\R)\times H^2_{\omega^{-1}}(\R)$. Applying the same argument, we also get $\U_*(x)+\tau\cV(x)\succ (0,0)$ for all $x\in\R$.

\textbf{Step 2.} Let $\tau_0>0$ be fixed. We use the asymptotics 
\bqs
\omega^{-1}(x)\U_*(x)=\omega^{-1}(x)\left(-U_*'(x),-W_*'(x) \right)\underset{+\infty}{\sim}\left(\beta\gamma_*,-\dfrac{rb \beta \gamma_*}{d_v^{+\infty}(\nu_u^{+\infty})} \right)x
\eqs
as $x\rightarrow+\infty$ to get the existence of $N>0$ large enough such that for all $\tilde{\tau}\in[\tau-\tau_0,\tau+\tau_0]$, one has
\bqs
\omega^{-1}(x)\U_*(x)\succ \omega^{-1}(x) |\tilde{\tau}\cV(x)|, \quad x\geq N, 
\eqs
such that it also holds by positivity of $\omega$ that
\bqs
\U_*(x)\succ |\tilde{\tau}\cV(x)|, \quad x\geq N.
\eqs
Now we use the fact $|\tau \cV(x)| \prec \U_*(x)$ for all $x\in\R$ to get that there exist $0<\epsilon\leq \tau_0$ small enough such that for all $\hat{\tau}\in(\tau-\epsilon,\tau+\epsilon)$
\bqs
\U_*(x)\succ |\hat{\tau}\cV(x)|, \quad x\in[-N,N].
\eqs
Combing the two estimates, we have that
\bqs
\U_*(x)\succ |\hat{\tau}\cV(x)|, \quad x\geq -N,
\eqs
for any $\hat{\tau}\in(\tau-\epsilon,\tau+\epsilon)$.

\textbf{Step 3.} We now fix $\hat{\tau}\in(\tau-\epsilon,\tau+\epsilon)$, and we want to show that 
\bqs
\U_*(x)\succ |\hat{\tau}\cV(x)|, \quad x\leq -N.
\eqs
To do so, we rewrite the operator $\tscL$ as follows
\bqs
\tscL=\left( \begin{array}{cc} 1 & 0 \\
0 & \sigma \end{array}\right)\partial_{xx}+c_*\left( \begin{array}{cc} 1 & 0 \\
0 & 1 \end{array}\right)\partial_x+\left( \begin{array}{cc} 1-a-2U_*(x)+aW_*(x) & aU_*(x) \\
rb(1-W_*(x)) & r(-1-bU_*(x)+2W_*(x)) \end{array}\right),
\eqs
and we denote
\bqs
\tscM(x):=\left( \begin{array}{cc} 1-a-2U_*(x)+aW_*(x) & aU_*(x) \\
rb(1-W_*(x)) & r(-1-bU_*(x)+2W_*(x)) \end{array}\right).
\eqs
As $x\rightarrow-\infty$, we have that
\bqs
\tscM(x)\rightarrow \tscM_-:=\left( \begin{array}{cc} -1 & a \\
0 & r(1-b) \end{array}\right).
\eqs
Upon denoting $P_+=(1,1)$ we have that $\tscM_-P_+\prec (0,0)$, and as a consequence (increasing $N$ if necessary), we have that $\tscM(x)P_+\prec (0,0)$ for all $x\leq -N$ such that $\tscL P_+ \prec (0,0)$ holds also true for $x\leq -N$.

Suppose that $\U_*(x)\succeq |\hat{\tau}\cV(x)|$ for all $x\leq -N$, then applying the maximum principle we get that $\U_*(x)\succ |\hat{\tau}\cV(x)|$ for all $x\leq -N$ and the conclusion follows. So let assume that there is some $\xi\in(-\infty,-N)$ such that one of the components of the vector $\U_*(x)-\hat{\tau}\cV(x)$ takes a local negative minimum at this point, say the first component. We consider the function $\cW_\iota(x):=\U_*(x)-\hat{\tau}\cV(x)+\iota P_+=(w_1^{\iota}(x),w_2^\iota(x))$ for $\iota>0$. As both $\U_*(x)$ and $\cV(x)$ converge to $(0,0)$ as $x\rightarrow-\infty$ there is a sufficiently large $\iota>0$ such that $\cW_\iota(x)\succ (0,0)$ for all $x\leq -N$. We now decrease $\iota$ up to $\iota_0>0$ such that 
\bqs
w_1^{\iota_0}(\xi)=0, \quad \cW_{\iota_0}(x)\succ (0,0) \text{ for } x\in(-\infty,-N), x\neq \xi.
\eqs
For such a $\iota_0$, we have that
\bqs
\tscL \cW_{\iota_0} = \iota_0 \tscM(x)P_+\prec (0,0), \quad x\in(-\infty,-N),
\eqs
but on the other hand we have that
\bqs
\partial_{xx}w_1^{\iota_0}(\xi)+c_*\partial_x w_1^{\iota_0}(\xi)+(1-a-2U(x)+aW_*(x))w_1^{\iota_0}(\xi)+aU_*(x)w_2^{\iota_0}(\xi)\geq 0,
\eqs
as $w_1^{\iota_0}$ attains its minimum at $\xi$ and that $w_2^{\iota_0}(x)\geq0$ on $x\leq-N$ by construction. This is a contradiction and consequently we have that $\U_*(x)-\hat{\tau}\cV(x)\succeq (0,0)$ for all $x\leq-N$ from which we deduce $\U_*(x)\succ |\hat{\tau}\cV(x)|$ for all $x\leq -N$ by the maximum principle. A similar argument shows that $\U_*(x)+\hat{\tau}\cV(x)\succeq (0,0)$ for all $x\leq-N$. 

In summary, we have shown that $\U_*(x)\succ |\hat{\tau}\cV(x)|$ for all $x\in\R$ which implies that $\mathcal{S}$ is open. As a consequence, we must have $\mathcal{S}=\R$ which contradicts the boundedness of $\U_*$. Thus $0$ is not an eigenvalue of $\tscL$ on  $L^2_{\omega^{-1}} (\R)\times L^2_{\omega^{-1}}(\R)$, which obviously implies that $0$ is not an eigenvalue of $\tcL$ on  $L^2_{\omega^{-1}} (\R)\times L^2_{\omega^{-1}}(\R)$ and completes the proof of Lemma~\ref{lem:0eig}.

\subsection{Proof of Lemma~\ref{lem:noposeig}}

First, we introduce the conjugate operator
\bqs
\scL = \omega^{-1} \tscL \omega, \quad \mathcal{D}(\scL)=H^2(\R)\times H^2(\R),
\eqs
with
\begin{align*}
\scL&=\left( \begin{array}{cc} 1 & 0 \\
0 & \sigma \end{array}\right)\partial_{xx}+\left( \begin{array}{cc} c_*+2\frac{\omega'}{\omega} & 0 \\
0 & c_*+2\sigma\frac{\omega'}{\omega} \end{array}\right)\partial_x+\scM(x), \\
\scM(x)&=\left( \begin{array}{cc} 1-a+c_*\frac{\omega'}{\omega}+\frac{\omega''}{\omega}-2U_*(x)+aW_*(x) & aU_*(x) \\
rb(1-W_*(x)) & r(-1+c_*\frac{\omega'}{r\omega}+\frac{\sigma\omega''}{r\omega}-bU_*(x)+2W_*(x)) \end{array}\right).
\end{align*}

Let assume that $\lambda>0$ is an eigenvalue of $\scL$ defined on $L^2(\R)\times L^2(\R)$ with corresponding eigenfunction $\cV$, that is
\bqs
\scL \cV = \lambda \cV, \quad 0\neq \cV \in H^2(\R)\times H^2(\R).
\eqs
Without of generality, we assume that the first component of $\cV=(v_1,v_2)$ is positive at some point $x_0$, that is $v_1(x_0)>0$. We claim that there exists some $\tau>0$ such that
\bqs
\cV(x)\preceq \tau \omega^{-1}(x)\U_*(x), \quad x \in\R,
\eqs
where $\U_*(x)=(-U_*'(x),-W_*'(x))\succ (0,0)$ satisfies $\tscL \U_*=0$ or equivalently $\scL\left(\omega^{-1}\U_*\right)=0$. First, as $\omega^{-1}(x)\U_*(x)$ is unbounded near $x=\infty$, we can always find $A>0$ large enough such that 
\bqs
\cV(x)\prec \omega^{-1}(x)\U_*(x), \quad x \geq A.
\eqs
Furthermore, the fact that $\U_*\succ (0,0)$ ensures that one can find $\tau>0$ such that 
\bqs
\cV(x)\preceq \tau \omega^{-1}(x)\U_*(x), \quad x\in[-A,+\infty).
\eqs
We are going to show that that the above inequality also holds true for $x\leq-A$. First, we remark that the matrix $\scM(x)$ asymptotically converges at $-\infty$ to the limiting matrix
\bqs
\scM_-=\left( \begin{array}{cc} -1+c_*\delta+\delta^2 & a \\
0 & r(1-b)+c_*\delta+\sigma\delta^2 \end{array}\right),
\eqs
and that $-1+c_*\delta+\delta^2<0$ and $r(1-b)+c_*\delta+\sigma\delta^2<0$ for our choice of $\delta$. As a consequence there exists a constant vector $Q_+\succ(0,0)$ such that $\scM_-Q_+\prec(0,0)$, which further implies that
\bqs
\scM(x)Q_+\prec(0,0), \quad x\leq -A,
\eqs
upon eventually increasing $A$. This also implies that one can find $\epsilon>0$, small enough such that
\bqs
\left(\left( \begin{array}{cc} 1 & 0 \\
0 & \sigma \end{array}\right)\epsilon^2-\left( \begin{array}{cc} c_*+\sigma\frac{\omega'}{\omega} & 0 \\
0 & c_*+2\sigma\frac{\omega'}{\omega} \end{array}\right)\epsilon+\scM(x) \right)Q_+ \prec(0,0), \quad x\leq -A.
\eqs
Assume now by contradiction that there exists some $x_1 < -A$, such that one of the components of $\tau \omega^{-1}(x)\U_*(x)-\cV(x)$ takes a local negative minimum at this point, say the first component. We introduce the function $\Q(x)=e^{-\epsilon x}Q_+\succ(0,0)$ which is unbounded as $x\rightarrow-\infty$, and further satisfies
\bqs
\scL \Q(x) = \left(\left( \begin{array}{cc} 1 & 0 \\
0 & \sigma \end{array}\right)\epsilon^2-\left( \begin{array}{cc} c_*+\sigma\frac{\omega'}{\omega} & 0 \\
0 & c_*+2\sigma\frac{\omega'}{\omega} \end{array}\right)\epsilon+\scM(x) \right)e^{-\epsilon x}Q_+\prec(0,0), \quad x\leq -A.
\eqs
As a consequence, one can find $\theta$ large enough such that
\bqs
\cW_\theta(x)=(w_1^\theta(x),w_2^\theta(x)):=\tau \omega^{-1}(x)\U_*(x)+\theta \Q(x)-\cV(x)\succ(0,0), \quad x\leq -A.
\eqs
We now decrease $\theta$ up to $\theta_0$ such that
\bqs
w_1^{\theta_0}(x_1)=0, \quad \cW_{\theta_0}(x)\succ(0,0), \quad x\leq -A, \quad x \neq x_1.
\eqs
On the one hand, we have that
\bqs
\scL \cW_{\theta_0} -\lambda \cW_{\theta_0}=-\tau \lambda \omega^{-1}\U_*+  \theta_0 \scL \Q- \lambda \theta_0 \Q \prec(0,0), \quad x\leq -A.
\eqs
On the other hand, the first component of $\scL \cW_{\theta_0} -\lambda \cW_{\theta_0}$ satisfies 
\begin{align*}
\left(\scL \cW_{\theta_0} -\lambda \cW_{\theta_0}\right)_1(x_1)=&~\partial_{xx} w_1^{\theta_0}(x_1)+\left(c_*+\sigma\frac{\omega'}{\omega}(x_1)\right)\partial_x w_1^{\theta_0}(x_1) \\
&~+\left(1-a-\lambda+c_*\frac{\omega'}{\omega}(x_1)+\frac{\omega''}{\omega}(x_1)-2U_*(x_1)+aW_*(x_1) \right)w_1^{\theta_0}(x_1)\\
&~+aU_*(x_1)w_2^{\theta_0}(x_1)\\
=&~\partial_{xx} w_1^{\theta_0}(x_1)+aU_*(x_1)w_2^{\theta_0}(x_1)\geq 0,
\end{align*}
which establishes a contradiction.

As a consequence, we have just proved that the set 
\bqs
\mathcal{E}=\left\{ \tau \geq 0 ~|~\cV(x)\preceq \tau \omega^{-1}(x)\U_*(x), \quad x \in\R \right\}
\eqs
is non-empty. We denote $\tau_0$ the greatest lower bound of $\mathcal{E}$. Suppose that $\tau_0>0$, then we have
\bqs
\cV(x)\preceq \tau_0 \omega^{-1}(x)\U_*(x), \quad x \in\R.
\eqs
Using the maximum principle as in Lemma~\ref{lem:0eig}, we get that necessarily 
\bqs
\cV(x)\prec \tau_0 \omega^{-1}(x)\U_*(x), \quad x \in\R.
\eqs
Now, repeating the argument that we just did to prove that $\mathcal{E}\neq \emptyset$, one can find $\delta>0$ small enough such that
\bqs
\cV(x)\preceq (\tau_0-\delta) \omega^{-1}(x)\U_*(x), \quad x \in\R,
\eqs
showing that $\tau_0-\delta \in \mathcal{E}$ contradicting the minimality of $\tau_0$. Hence, we must have $\tau_0=0$, which then implies that $\cV(x)\preceq(0,0)$ contradicting the fact $v_1(x_0)>0$. Thus we must have $\lambda \leq 0$, but from Lemma~\ref{lem:0eig}, $\lambda=0$ is not an eigenvalue of $\scL$, such that $\lambda<0$. This concludes the proof of the lemma.

\bibliographystyle{abbrv}
\bibliography{LVBib}

\end{document}